\documentclass[12pt,leqno]{amsart}  
\usepackage{pgf,pgfarrows,pgfshade}		
\usepackage{amsmath,amstext,amsthm,amssymb,amsxtra}
\usepackage[colorlinks,pagebackref,hypertexnames=false]{hyperref}
\usepackage[backrefs]{amsrefs}
\allowdisplaybreaks
\swapnumbers
   
\theoremstyle{plain} 
\newtheorem{lemma}[equation]{Lemma} 
\newtheorem{proposition}[equation]{Proposition} 
\newtheorem{theorem}[equation]{Theorem}

\theoremstyle{definition}

\theoremstyle{remark}
\newtheorem{remark}[equation]{Remark}
\newtheorem*{Acknowledgment}{Acknowledgment}

\def\thm#1/{Theorem~\ref{t#1}}
\def\c#1/{Corollary~\ref{c#1}}
\def\l#1/{Lemma~\ref{l#1}}
\def\s#1/{Section~\ref{s#1}}
\def\e#1/{(\ref{e#1})}
\def\d#1/{Definition~\ref{d#1}}

\def\Label #1 {\label{#1}}

\def\norm#1.#2.{\lVert#1\rVert_{#2}}
\def\Norm#1.#2.{\bigl\lVert#1\bigr\rVert_{#2}}
\def\NOrm#1.#2.{\Bigl\lVert#1\Bigr\rVert_{#2}}
\def\NORm#1.#2.{\biggl\lVert#1\biggr\rVert_{#2}}
\def\NORM#1.#2.{\Bigggl\lVert#1\Bigggr\rVert_{#2}}

\def\ip#1,#2.{\langle #1,#2\rangle}
\def\Ip#1,#2.{\bigl\langle#1,#2\bigr\rangle}
\def\IP#1,#2.{\Biggl\langle#1,#2\Biggr\rangle}

\def\abs#1{\lvert#1\rvert}
\def\Abs#1{\bigl\lvert#1\bigr\rvert}
\def\ABs#1{\biggl\lvert#1\biggr\rvert}

\def\Xint#1{\mathchoice
   {\XXint\displaystyle\textstyle{#1}}%
   {\XXint\textstyle\scriptstyle{#1}}%
   {\XXint\scriptstyle\scriptscriptstyle{#1}}%
   {\XXint\scriptscriptstyle\scriptscriptstyle{#1}}%
   \!\int}
\def\XXint#1#2#3{{\setbox0=\hbox{$#1{#2#3}{\int}$}
     \vcenter{\hbox{$#2#3$}}\kern-.5\wd0}}

\def\dashint{\Xint-}

 \def\ind#1{{\mathbf 1}_{#1}}

 \def\Mod#1{\operatorname{Mod}_{#1}}
 \def\Tr#1{\operatorname{Tr}_{#1}}
 \def\Dil#1#2{\operatorname{Dil}_{#2}^{(#1)}}
 \def\Well#1{\operatorname{Well}(#1)}
 \def\sh#1{\operatorname{sh}(#1) }
 \def\emb#1 {\operatorname{emb}(#1)} 
 \def\Enl#1  {\operatorname{Enl} (#1)}

\def\mid{\,:\,}

\begin{document}

\title[Rubio de Francia's Littlewood--Paley Inequality]
{Issues related to Rubio de Francia's \\
Littlewood--Paley Inequality}
\author[M.T. Lacey]{Michael T. Lacey\\
School of Mathematics\\
Georgia Institute of Technology\\
Atlanta GA 30332\\
}

\thanks{Research supported in part by a National 
Science Foundation Grant. The author is a Guggenheim Fellow.} 

\email{lacey@math.gatech.edu}

\begin{abstract}  Let $\operatorname S_\omega f=\int_\omega \widehat f(\xi)e^{ix\xi}\; d\xi$ be the Fourier projection 
operator to an interval $\omega$ in the real line.  Rubio de Francia's Littlewood Paley inequality \cite{rubio} states that 
for any collection of disjoint intervals $\Omega$, we have 
 \begin{equation*}
\NOrm \Bigl[ \sum_{\omega\in\Omega} \abs{\operatorname S_\omega f}^2\Bigr]^{1/2} .p.\lesssim{}\norm f.p.,\qquad 2\le{}p<\infty.
 \end{equation*}
We survey developments related to this inequality, including the higher dimensional case, and consequences for multipliers.
\end{abstract}

\maketitle
\tableofcontents

 \section{Introduction}

 Our subject is a group of topics related to Rubio de Francia's extension \cite{rubio} of the classical 
 Littlewood Paley inequality.  We are especially interested in presenting a proof that 
 highlights an approach in the language of time--frequency analysis, and  
 addresses the known higher dimensional versions of this Theorem.  It is hoped that this 
 approach will be helpful in conceiving of new versions of these inequalities. 
 A first result in this direction is in the result of Karagulyan and the author 
 \cite {math.CA/0404028}. 
  These inequalities yield interesting consequence for multipliers, and these are reviewed as well.

 Define the Fourier transform by 
  \begin{equation*}
 \widehat f (\xi)=\int_{\mathbb R^d} f(x)e^{-ix\cdot\xi}\; dx.
  \end{equation*} 
 In one dimension, the projection onto the positive frequencies 
  \begin{equation*}
\operatorname  P_+f(x):=\int_0^\infty \widehat f(\xi) e^{ix \xi} \; d\xi
  \end{equation*}
 is a bounded operator on all $L^p(\mathbb R)$, $1<p<\infty$. 
 The typical proof of this fact first establishes the $L^p$ inequalities for the Hilbert transform, given by, 
  \begin{equation*}
 \operatorname  Hf(x):=\lim_{\epsilon\to\infty} \int_{\abs y>\epsilon} f(x-y)\frac {dy}y
  \end{equation*} 
 The Hilbert transform is given in frequency by a constant times 
  \begin{equation*}
\operatorname h f(x)=c \int \widehat f(\xi) \,\text{sign}(\xi)\operatorname e^{ix\xi}\; d\xi.
  \end{equation*}
 We see that $\operatorname P_+$ is linear combination of the identity 
 and $\operatorname H$.   In particular $\operatorname P_+$ and $\operatorname H$ enjoy the same mapping properties.

 In this paper, we will take the view that 
 $L^p(\mathbb R^d)$ is the tensor product of $d$ copies of $L^p(\mathbb R)$.  
 A particular consequence is that the 
 projection onto the positive quadrant 
  \begin{equation*}
 \operatorname P_+f(x):=\int_{[0,\infty]^d} f(\xi) \operatorname  e^{ix \cdot\xi} \; d\xi
   \end{equation*}
  is a bounded operator on all $L^p(\mathbb R^d)$, as it is merely a tensor product of the one dimensional projections.

     A 
  rectangle in $\mathbb R^d$ is denoted by $\omega$.  Define the Fourier restriction operator to be 
   \begin{equation*}
  \operatorname S_\omega f(x)=\int_\omega \operatorname  e^{ix\cdot\xi} \widehat f(\xi)\; d\xi,
   \end{equation*}
 This projection operator is bounded on all $L^p(\mathbb R^d)$, with constant bounded 
 independently of $\omega$.  To see this, define the 
 modulation operators by 
  \begin{equation}\label{e.moddef}
 \Mod \xi f(x):= \operatorname e^{ix\cdot \xi} f (x)
  \end{equation}
 Observe that for $\xi=(\xi_1,\ldots,\xi_d)$, the interval
 $\omega=\prod_{j=1}^d [\xi_j,\infty)$, we have 
 $\operatorname S_\omega=\Mod {-\xi} \operatorname P_+ \Mod {\xi}$. 
 Hence this projection is uniformly bounded.  By taking linear combinations of 
 projections of this type, we can obtain the $ L^p$ boundedness of any 
 projection operator $ \operatorname S _{\omega }$, for rectangles $ \omega $.
 
 The Theorem we wish to explain is

   \begin{theorem}\label{t.lp}  Let $\Omega$ be any collection of disjoint rectangles
    with respect to a fixed choice of basis. Then the square
  function below maps $L^p(\mathbb R^d)$ into itself for $2\le p<\infty$.
   \begin{equation*}
  \operatorname S^\Omega f(x):=\Bigl[\sum_{\omega\in\Omega}\abs{\operatorname S_\omega f(x)}^2\Bigr]^{1/2}.
   \end{equation*}
   \end{theorem}
  
  In one dimension this is Rubio de Francia's Theorem \cite{rubio}.  His proof 
  pointed to the primacy of a $BMO$ estimate in the proof 
  of the Theorem.
  The higher dimensional form was investigated by J.-L.~Journ\'e \cite{MR88d:42028}.  His original argument has been 
  reshaped by F.~Soria \cite{MR88g:42026},   S.~Sato, \cite{MR92c:42020}, and Xue~Zhu \cite{MR93f:42041}.  In this instance, the product 
  $BMO$ is essential, in the theory as developed by S.-Y.~Chang and R.~Fefferman 
  \cites{MR81c:32016,MR82a:32009,MR86g:42038}.   

   We begin our discussion with the one dimensional case, 
   followed by  the higher dimensional case.  We adopt a `time-frequency' approach to the 
   Theorem, inspired in part by the author's joint work with Christoph Thiele 
   \cites{laceythielecarleson
,
MR99b:42014}..
   The same pattern is adopted for the 
   multiplier questions.  The paper concludes with notes and comments.

   We do not keep track of the value of generic absolute constants, instead
   using the  notation $A\lesssim{}B$ iff $A\le{}KB$ for some constant $K$.  
   Write $A\simeq B$ iff $A\lesssim{}B$ and $B\lesssim{}A$. For a rectangle 
   $\omega$ and scalar $\lambda>0$, $\lambda\omega$ denotes the rectangle with the 
   same center as $\omega$ but each side length is $\lambda$ times the same side length of $\omega$.   We use the notation
$\ind A$ to denote the indicator function of the set $A$, that is, $ 1_A(x)=1$ if 
$ x\in A$ and is otherwise $ 0$. 
Averages of integrals over a set are written as 
 \begin{equation*}
\dashint_A f\;dx:=\abs{A}^{-1}\int_A f\;dx.
 \end{equation*}
For an operator $T$, $\norm T.p.$ denotes the norm of $T$ as an operator from $L^p(\mathbb R^d)$ to itself.  In addition to the Modulation operator defined above, we will also use the translation operator 
 \begin{equation*}
\Tr y f(x):=f(x-y).
 \end{equation*}
We shall assume the reader is familiar with the norm bounds for the one dimensional maximal function 
 \begin{equation*}
\operatorname Mf(x)=\sup_t\dashint_{[-t,t]} \abs{f(x-y)}\; dt
 \end{equation*}
The principal fact we need is that it  maps $L^p$ into itself for $1<p<\infty$.  In $d$ dimensions, the strong maximal function refers to the 
maximal function 
 \begin{equation*}
\operatorname Mf(x)=\sup_{t_1,\ldots,t_d>0}\dashint_{[-t_1,t_1]\times\cdots[-t_d,t_d]} \abs{f(x_1-y_1,\ldots,x_d-y_d)} \; dy_1\cdots dy_d
 \end{equation*}
Note that this maximal function is less than the one dimensional maximal function applied in each coordinate in succession.

\begin{Acknowledgment}
An initial version of these notes was prepared while in residence at the Schr\"odinger 
Institute of Vienna Austria.  The paper has been improved by the efforts of a 
 conscientious  referee. 
\end{Acknowledgment}

   \section{The One Dimensional Argument}
   In this setting, we give the proof in one dimension, as it is very much easier in this case.  In addition, some of the ideas 
   in this case will extend immediately to the higher dimensional case.    
   
   \subsection{Classical Theory}

   We should take some care to recall the classical theory of Littlewood and Paley.  Let $\Delta$ denote the dyadic intervals 
    \begin{equation*}
   \Delta:=\{\epsilon[2^k,2^{k+1})\mid \epsilon\in\{\pm1\},\ k\in\mathbb Z\}.
    \end{equation*}
   
   The classical Theorem is that 
   
    \begin{theorem}\label{t.lpineq} 
   For all $1<p<\infty$, we have 
    \begin{equation} \label{e.lp}
   \Norm \operatorname S^\Delta f.p.\simeq\norm f.p.
    \end{equation}
    \end{theorem}
   
   We will not prove this here, but will make comments about the proof.  If one knows that 
    \begin{equation} \label{e.upper}
   \Norm \operatorname S^\Delta f.p.\lesssim\norm f.p., \qquad 1<p<\infty
    \end{equation}
   then a duality argument permits one to deduce the reverse inequality for $L^{p'}$ norms,  $p'=p/(p-1)$.  Indeed, 
   for $g\in L^{p'}$, choose $f\in L^{p}$  of norm one so that $\norm g.p'.=\ip f,g.$.  Then 
    \begin{equation} \label{e.duality}
\begin{split}
   \norm g.p'.={}&\ip f,g.
   \\{}={}&\int 
   \sum_{\omega\in\Delta}  \operatorname S_\omega f\overline{\operatorname S_\omega g}\; dx \nonumber
   \\{}\le{}& \ip \operatorname S^\Delta f, \operatorname S^\Delta g.    \nonumber
   \\{}\le{}& \norm \operatorname S^\Delta f.p. \norm \operatorname S^\Delta g.p'.   \nonumber
   \\{}\lesssim{}& \norm \operatorname S^\Delta g.p'.   \nonumber
\end{split}\end{equation}
  One only need prove the upper inequality for the full range of $1<p<\infty$.

  In so doing, we are faced with a common problem in the subject. 
  Sharp frequency jumps produce  kernels with slow decay at infinity,
  as is evidenced by the Hilbert transform, which has a single frequency jump
   and a non--integrable kernel.  The operator $\operatorname S^\Delta$ has infinitely
   many frequency jumps. It is far easier to  to 
  study a related operators with smoother frequency behavior, for then  standard aspects of  
  Calder\'on--Zygmund Theory at one's disposal.
  Our purpose is then to introduce a class of operators which mimic the behavior of $\operatorname S^\Delta$, but 
  have smoother frequency behavior. 
  
   Consider a smooth function $\psi_+$ which satisfies $\ind {[1,2]}\le\widehat {\psi_+}\le\ind {[\tfrac12,\tfrac52]}$. 
   Notice that $\psi*f$ is a smooth version of $S_{[1,2]}f$.  
   Let $\psi_-=\overline{\psi_+}$. Define the dilation operators, of scale $ \lambda $,  by 
    \begin{equation}\label{e.dilate}
   \Dil p\lambda f(x):= \lambda^{-1/p}f(x/\lambda),\qquad 0<p\le\infty,\ \lambda>0.
    \end{equation}
  The normalization chosen here normalizes the $ L^p$ norm of $ \Dil p \lambda $  to be one.

   Consider distributions of the form 
    \begin{equation} \label{e.K}
   K=\sum_{k\in\mathbb Z}\sum_{\sigma\in\{\pm\}} \varepsilon_{k,\sigma}\Dil 1{2^k}\psi_\sigma,\qquad \varepsilon_{k,\sigma}\in\{\pm1\}.
    \end{equation}
   and the operators $\operatorname Tf=K*f$.  This class of distributions satisfy the standard estimates of Calder\'on--Zygmund theory, with 
   constants independent of the choices of signs above. 
   In particular, these estimates would be 
    \begin{gather*}
   \sup_\xi \abs{\widehat K(\xi)}<C\,,
   \\
   \abs{K(y)}<C\abs{y}^{-1}\,,
   \\
   \abs{ \tfrac d{dy} K(y)}<C \abs{ y } ^{-2}\,,
    \end{gather*}
   for a universal constant $C$.  These inequalities imply that 
    the operator norms of $\operatorname T$  on $ L^p$ are bounded by constants that
    depend only on $p$. 
    
    The uniformity of the constants in the operator norms permits us to 
   average over the choice of signs, and apply the Khintchine inequalities to conclude that 
    \begin{equation} \label{e.smoothlp}
   \NOrm \Bigl[\sum_{k\in\mathbb Z}\sum_{\sigma\in\{\pm\}} \abs{\Dil 1{2^k}\psi_\sigma *f }^2 \Bigr]^{1/2} .p.\lesssim{}\norm f.p.,\qquad 1<p<\infty.
    \end{equation}
   This is nearly the upper half of the inequalities in \thm.lp/.  For historical reasons, ``smooth'' square functions such as the one above,
   are referred to as ``$G$ functions.''
   
To conclude the Theorem as stated, one method uses an extension of the boundedness of the  Hilbert transform to a vector valued setting. 
The particular form needed concerns the extension of the Hilbert transform to functions taking values in $\ell^q$ spaces.  In particular, we have the inequalities 
 \begin{equation} \label{e.vectorH}
\lVert \norm  {\operatorname H f_k} .\ell^q . \rVert_p\lesssim{}C_{p,q}\lVert \norm f_k.\ell^q. \rVert_p,\qquad 1<p,q<\infty.
 \end{equation}
Vector valued inequalities are strongly linked to weighted inequalities, and one of the standard approaches to these inequalities depends upon 
the beautiful inequality of C.~Fefferman and E.M.~Stein \cite{MR44:2026}
 \begin{equation}\label{e.weighted}
\int \abs{\operatorname Hf}^q g\; dx\lesssim{} \int \abs{f}^q 
(\operatorname M\abs{g}^{1+\epsilon})^{1/(1+\epsilon)}\; dx, 
\qquad 1<q<\infty, \ 0<\epsilon<1.
 \end{equation}
The implied constant depends only on $q$ and $\epsilon$.  While we stated this for the Hilbert transform, it is important for our purposes 
to further note that this inequality continues to hold for a wide range of Calder\'on--Zygmund operators, including those that occur in (\ref{e.K}).
This is an observation that goes back to J.~Schwartz \cite{MR0143031}, with 
many extensions, especially that of  Benedek, Calder{\'o}n and Panzone \cite{MR0133653}.

The proof that (\ref{e.weighted}) implies (\ref{e.vectorH}) follows.  Note that we need only prove the vector valued estimates for 
$1<q\le{}p<\infty$, as the  remaining estimates follow by duality, namely  the dual estimate of 
$\operatorname H\mid L^p(\ell^q)\longrightarrow L^p(\ell^q)$ is 
$\operatorname H\mid L^{p'}(\ell^{q'})\longrightarrow  L^{p'}(\ell^{q'})$, 
in which the primes denote the conjugate index, $ p'=p/(p-1)$.
 The cases of $q=p$ are trivial. For $1<q<{}p<\infty$, and  $ \{f_k\}\in L^p(\ell^q)$ of norm one, it suffices to show that 
  \begin{equation*}
 \NOrm \sum_k \abs{\operatorname H f_k}^q .p/q.\lesssim{}1.
  \end{equation*}
 To do so, by duality, we can take $g\in L^{(p/q)'}$ of norm one, and estimate 
  \begin{align*}
 \sum_k\int  \abs{\operatorname H f_k}^qg\; dx\lesssim{}& \sum_k \int \abs {f_k}^q (\operatorname M \abs{g}^{1+\epsilon})^{1/(1+\epsilon)} \; dx
 \\{}\lesssim{}& \NOrm \sum_k \abs {f_k}^q .p/q.\norm (\operatorname M \abs{g}^{1+\epsilon})^{1/(1+\epsilon)} .(p/q)'.
 \\{}\lesssim{}& 1
  \end{align*}
 provided we take $1+\epsilon< (p/q)'$.

Now, the Fourier  projection onto  an interval $ \omega $ can be 
obtained as a linear combination of modulations of the Hilbert transform. 
Using this, one sees that  the estimate (\ref{e.vectorH}) 
extends to the Fourier projections onto intervals. 
Namely, we have the estimate 
 \begin{equation*}
\lVert \norm \operatorname S_\omega f_\omega .\ell^2(\Omega). \rVert_p\lesssim{} \lVert \norm f_\omega. \ell^2(\Omega). \rVert_p,\qquad 1<p<\infty. 
 \end{equation*}
This is valid for all collections of intervals $\Omega$.  Applying it to (\ref{e.smoothlp}), with $\Omega=\Delta$, and using the fact that 
$S_{\sigma[2^k,2^{k+1})}f= \operatorname S_{\sigma[2^k,2^{k+1})}\Dil 1{2^k} \psi_\sigma*f$ proves the upper half of the inequalities of \thm.lpineq/, 
which what we wanted.

\medbreak 

For our subsequent use, we note that the vector valued extension of the Hilbert transform depends upon structural estimates that continue to hold for a wide variety of Calder\'on--Zygmund kernels.  In particular, the Littlewood--Paley inequalities also admit a vector valued extension, 
 \begin{equation}\label{e.lpvector}
\lVert \norm \operatorname S^\Delta f_k .\ell^q.\rVert_p\simeq \lVert \norm f_k.\ell^q. \rVert_p,\qquad 1<p,q<\infty. 
 \end{equation}

\subsection{Well--Distributed Collections}

We begin the main line of argument for Rubio de Francia's inequality in one dimension.  
The first step, found by Rubio de Francia \cite {rubio}, is a reduction of the general case 
to one in which one can square function by a smoother object.

Say that a collection of intervals $\Omega$ is {\em well distributed} if 
 \begin{equation}\label{e.well}
\NOrm \sum_{\omega\in\Omega}\ind {3\omega}.\infty.\le100.
 \end{equation}
 Thus, after dilating the intervals in the collection of a factor of (say) $ 3$, 
 at most $100$ intervals can intersect. 
 
The well distributed collections allow one to smooth out $\operatorname S_\omega$,
just as one does $\operatorname S_{[1,2]}$ in the proof of the classical 
Littlewood--Paley inequality. The main fact we should observe here is that 

 \begin{lemma}\label{l.wellsuffice}  For each collection of intervals $\Omega$, we can define a well distributed collection $\Well \Omega$ for which 
 \begin{equation*}
\norm \operatorname S^\Omega f.p.\simeq \norm \operatorname S^{\Well \Omega}f.p.,\qquad 1<p<\infty. 
 \end{equation*}
 \end{lemma}

\begin{proof}

The argument here depends upon inequalities for vector valued singular integral operators. 
We define the collection $\Well \Omega$ by first considering the  interval  $[-\frac12,\frac12]$. Set
 \begin{equation*}
\Well {[-\tfrac12,\tfrac12]}=\{ [-\tfrac1{18},\tfrac1{18}] ,\pm[\tfrac12-\tfrac4{9}(\tfrac45)^{k},\tfrac12-\tfrac4{9}(\tfrac45)^{k+1}]\mid k\ge0\}.
 \end{equation*}
It is straightforward to check that all the  intervals in this collection have a distance to the boundary of $[-\frac12,\frac12]$ that is 
four times their length. 
In particular, 
this collection is well distributed.  It has the additional property that 
for each $\omega\in \Well {[-\frac12,\frac12]}$ we have 
$2\omega\subset [-\frac12,\frac12]$.    

  It is an extension of the usual Littlewood--Paley inequality that  
 \begin{equation*}
\norm \operatorname S_{[-1/2,1/2]}f.p.\simeq \norm \operatorname S^{\Well {[-1/2,1/2]}} S_{[-1/2,1/2]}f.p.,\qquad 1<p<\infty.
 \end{equation*}
This inequality continues to hold in the  vector valued setting of (\ref{e.lpvector}). 

We define $\Well \omega $ by affine invariance. For an interval $\omega$, select an affine function $\alpha\mid [-\frac12,\frac12]\longrightarrow \omega$, we set 
$\Well \omega:=\alpha(\Well {[-\frac12,\frac12]})$.  
For collections of intervals $ \Omega $,  we define $\Well \Omega:=\bigcup_{\omega\in\Omega}\Well \omega$. It is clear that $\Well \Omega$ is well distributed for collections of disjoint intervals $\Omega$. By a vector valued Littlewood--Paley inequality, we have 
 \begin{equation*}
\norm \operatorname S^\Omega f.p.\simeq \norm \operatorname S^{\Well \Omega}f.p.,\qquad 1<p<\infty.
 \end{equation*}
This completes the proof of our Lemma.

\end{proof}

In the proof of the Lemma, we see that we are `resolving the frequency jump' at both 
endpoints of the interval.  In the sequel however, we don't need to rely upon this 
construction, using only the general definition of well--distributed.

For the remainder of the proof, we assume that $ \Omega $ is well--distributed.  
We need only consider a smooth version of the square function $ \operatorname S ^{\Omega }$, 
with the assumption of well distributed is  critical to boundedness of 
the smooth operator on $L^2$.

Let $\varphi$ be a Schwartz function 
so that 
 \begin{equation}\label{e.zvf}
\ind {[-1/2,1/2]}\le\widehat\varphi\le\ind {[-1,1]}
 \end{equation}
Set $\varphi^\omega=\Mod {c(\omega)}\Dil 2{\abs{\omega}^{-1}}\varphi$, and  
 \begin{equation*}
\operatorname G^\Omega f=\Bigl[ \sum_{\omega\in\Omega}\abs{\varphi^\omega*f}^2\Bigr]^{1/2}.
 \end{equation*}
We need only show that 
  \begin{equation} \label{e.G} 
\norm \operatorname G^\Omega f.p.\lesssim{}\norm f.p.,\qquad 2\le p<\infty,
 \end{equation}
for well distributed collections $\Omega$.  Note that that the well distributed assumption   and the assumptions about $\varphi$ make the $L^2$ inequality obvious. 

\subsection{The Tile Operator}

We use the previous Lemma to pass to an operator that is easier to control than 
the projections $\operatorname S_\omega$ or $\varphi^\omega*f$.
This is done  in the 
time frequency plane.  Let $ \mathbf D$ be the dyadic intervals in $\mathbb R$, that is 
 \begin{equation*}
\mathbf D:=\{[j2^k,(j+1)2^k\mid j,k\in\mathbb Z\}.
 \end{equation*}

Say that $s=I_s\times \omega_s$ is a {\em tile} if $I_s\in\mathbf D$, $\omega_s$ is an interval, and $1\le\abs s=\abs{I_s}\cdot\abs{\omega_s}<2$.  Note that for any $\omega_s$, there is 
one choice of $\abs{I_s}$ for which $I_s\times \omega_s$ will be a tile. 
We fix a Schwartz function $\varphi$, and define 
 \begin{equation*}
\varphi_s:=\Mod {c(\omega_s)} \Tr {c(I_s)}\Dil 2{\abs{I_s}} \varphi
 \end{equation*}
where $c(J)$ denotes the center of $J$. We take $\varphi$ as above, a Schwartz function 
satisfying  
$\mathbf 1 _{[-1,1]}\le{} \widehat \varphi\le{} \mathbf 1 _{[-2,2]}$. 

Choosing tiles to have  area approximately equal to one is suggested by the 
Fourier uncertainty principle. We  sometimes refer to  $ I_s$ and 
$ \omega _s$ as \emph {dual} intervals.
With this choice of definitions, the function $\varphi_s$ is approximately localized
in the time frequency plane to the rectangle $I_s\times \omega_s$. 
This localization is precise in the frequency variable.  
The function $\widehat {\varphi_s}$ is supported in the interval $2\omega_s$. 
But, $\varphi_s$ is only approximately 
supported near the interval $I_s$.  Since $\varphi$ is rapidly decreasing, we trivially have the estimate 
 \begin{equation*}
\abs{\varphi_s(x)}\lesssim{} \abs{I_s}^{-1/2} (1+\abs{I_s}^{-1}\,\abs{x-c(I_s)})^{-N},
\qquad N\ge1.
 \end{equation*}
This is an adequate substitute for being compactly supported in the time variable.

For a collection of intervals $\Omega$, we set $\mathcal T(\Omega)$ to be the set 
of all possible tiles $s$ such that $\omega_s\in\Omega$. 
Note that for each $\omega\in\Omega$, the set of intervals 
$\mathcal T( \{ \omega \})=\{I\mid I\times \omega\in\mathcal T(\Omega)\}$ is a a partition of $\mathbb R$ into intervals of equal length.   See Figure $1$. 
Associated to $\mathcal T(\Omega)$ is a natural square function 
 \begin{equation*}
\operatorname T^\Omega f=\Bigl[\sum_{s\in\mathcal T(\Omega)} \frac{\abs{\ip f,\varphi_s.}^2}{{\abs{I_s}}} \ind {I_s}\Bigr]^{1/2}.
 \end{equation*}

\begin{figure}
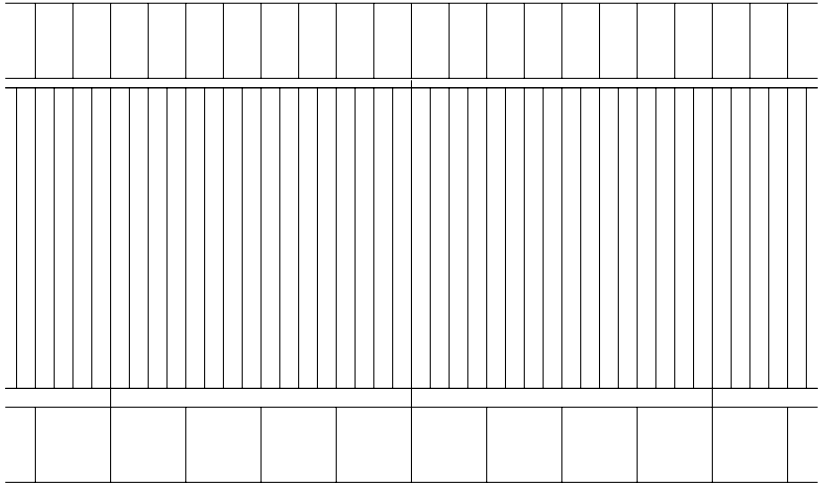

\begin{center}
\begin{pgfpicture}{0cm}{0cm}{8cm}{6cm}
\begin{pgftranslate}{\pgfxy(4,5.375)}  
\pgfgrid[stepx=0.5cm,stepy=1cm]{\pgfxy(-5.4,0)}{\pgfxy(5.4,1)}
\end{pgftranslate}
%
\begin{pgftranslate}{\pgfxy(4,5.25)}
\pgfgrid[stepx=8cm,stepy=0.125cm]{\pgfxy(-5.4,0)}{\pgfxy(5.4,0.1)}
\end{pgftranslate}
\begin{pgftranslate}{\pgfxy(4,1.25)}
\pgfgrid[stepx=.25cm,stepy=4cm]{\pgfxy(-5.4,0)}{\pgfxy(5.4,4)}
\end{pgftranslate}
%
\begin{pgftranslate}{\pgfxy(4,0)}
\pgfgrid[stepx=4cm,stepy=.25cm]{\pgfxy(-5.4,1)}{\pgfxy(5.4,1.25)}
\pgfgrid[stepx=1cm,stepy=1cm]{\pgfxy(-5.4,0)}{\pgfxy(5.4,1)}
\end{pgftranslate}
\end{pgfpicture}
\end{center}
\caption{Example tiles in a collection $\mathcal T(\Omega)$.  The vertical axis represents 
frequency, the horizontal is time. }
\end{figure}


Our main Lemma is that 

 \begin{lemma}\label{l.main1} 
For any collection of well distributed intervals $\Omega$, we have 
 \begin{equation*}
\norm \operatorname T^\Omega f.p.\lesssim{}\norm f.p.,\qquad 2\le p<\infty.
 \end{equation*}
 \end{lemma}

Let us argue that this Lemma proves (\ref{e.G}), for a slightly different square function, and so proves Rubio de Francia's 
Theorem in the one dimensional case.  One task is to pass 
from a sum of rank one operators to a convolution operator. This is in fact 
a general principle, that we can formulate this way.

 \begin{lemma}\label{l.convolve}  Let $\varphi$ and $\phi$ be real valued Schwartz functions on $\mathbb R$.  Then, 
 \begin{align*}
\dashint_{[0,1]}\sum_{m\in\mathbb Z} 
\ip f, \Tr {y+m}  \varphi. \Tr {y+m}  \phi \; dy=f*\Phi
\\
\text{where}\qquad \Phi(x)=\int \overline{\varphi(u)}\phi(x+u)\; du.
 \end{align*}
In particular, $\widehat{\Phi}=\overline{\widehat\varphi}\widehat \phi$. 
 \end{lemma}

The proof is immediate.  The integral in question is 
 \begin{equation*}
\iint_{\mathbb R} f(z)\overline{\varphi(z-y)}\phi(x-y)\; dydz
 \end{equation*}
and one changes variables, $u=z-y$.

\medskip

\begin{proof}[Proof of (\ref{e.G}).]
We need to pass from the discrete operator 
to a square function of convolution operators. 
Let 
 \begin{equation*}
\chi(x):=(1+\abs x)^{-10}, \qquad \chi_{(I)}=\Dil 1 {\abs I}\Tr {c(I)} \chi, 
 \end{equation*}
and set for $\omega\in\Omega$, 
 \begin{equation*}
\operatorname H_\omega f=\sum_{\substack{s\in\mathcal T(\Omega)\\ \omega_s=\omega}}
\ip f,\varphi_s. \varphi_s \,.
 \end{equation*}
By Cauchy--Schwarz, we may dominate 
 \begin{align*}
\abs{\operatorname H_\omega f}\le{}&
\sum_{\substack{s\in\mathcal T(\Omega)\\ \omega_s=\omega} }\abs{\ip f,\varphi_s. \varphi_s }
\\
{}\lesssim{}& \sum_{\substack{s\in\mathcal T(\Omega)\\ \omega_s=\omega}} \frac{\abs{\ip f,\varphi_s.}}{\sqrt{\abs{I_s}}} 
\abs{ \chi_{(I_s)}\ast \ind {I_s} }^2
\\ 
{}\lesssim{}& 
\Bigl[ \sum_{\substack{s\in\mathcal T(\Omega)\\ \omega_s=\omega} }
\frac{\abs{\ip f,\varphi_s.}^2}{{\abs{I_s}}} 
\abs{ \chi_{(I_s)}\ast \ind {I_s} }
\Bigr]^{1/2}.
 \end{align*}
We took some care to include the convolution in this inequality, so that we could use the easily verified 
inequality $\int \abs{\chi_{(I)}\ast f}^2g\; dx\le{}\int \abs{f}^2\chi_{(I)}\ast g\; dx$ in the
following way:  
The square function $\norm\operatorname H_\omega  f.\ell^2(\Omega).$ is seen to map $L^p$ into itself, $2<p<\infty$ by duality. For 
functions $g\in L^{(p/2)'}$ of norm one, we can estimate 
 \begin{align*}
 \sum_{s\in\mathcal T(\Omega)} \frac{\abs{\ip f,\varphi_s.}^2}{{\abs{I_s}}} 
\int \abs{ \chi_{(I_s)}\ast \ind {I_s} }g\;dx\le{}& 
 \sum_{s\in\mathcal T(\Omega)} \frac{\abs{\ip f,\varphi_s.}^2}{{\abs{I_s}}} 
\int \ind {I_s}\chi_{(I_s)}\ast g\;dx
\\{}\le{}&    
\int \abs{ \operatorname T^\Omega f}^2 \sup_I \chi_{(I)}\ast g\;dx 
\\{}\lesssim{}&  \lVert  \operatorname T^\Omega f
\rVert_p^2\norm \operatorname Mg. (p/2)'. 
\\{}\lesssim{}&  \lVert f \rVert_p^2.
 \end{align*}
Here, $(p/2)'$ is the conjugate index to $p/2$, and $M$ is the maximal function.

Thus, we have verified that 
 \begin{equation*}
\lVert \norm\operatorname H_\omega  f.\ell^2(\Omega). \rVert_p\lesssim{}\norm f.p.,\qquad 2<p<\infty. 
 \end{equation*}
We now  derive a convolution inequality.  By \l.convolve/, 
 \begin{equation*}
\lim_{T\to\infty} \dashint_{[0,T]} \Tr {-y}\operatorname H_\omega  \Tr y f \; dy=\psi^\omega*f,
 \end{equation*}
for all $\omega$, where $\widehat {\psi^\omega}=\abs{ \widehat {\varphi^\omega}}^2$. 

\end{proof}

Thus, we 
see that a square function inequality much like that of (\ref{e.G}) holds; 
 this completes the proof 
of Rubio de Francia's Theorem in the one dimensional case, aside from the proof of \l.main1/.

\subsection{Proof of \l.main1/}

The proof of the boundedness of the tile operator $ \operatorname T ^{\Omega }$ on $ L^2$ 
is straight forward, yet finer facts about this boundedness are very useful 
in extending the boundedness to $ L^p$ for $ p>2$.   This is the subject of the 
next Proposition.

\begin{proposition}\label{p.L2} 
Let $ \psi $ be a smooth, rapidly decreasing function, satisfying in particular 
\begin{equation}\label{e.zxc}
\Abs{\operatorname  \psi (x) } \lesssim (1+\abs x) ^{-20}, 
\end{equation}
For any interval $ \omega $, we have 
\begin{align} \label{e.zxc0}
\sum _{s\in \mathcal T(\{\omega \})} \abs{ \ip f,  \psi _s .  } ^{2} &\lesssim \norm f.2.^2
\end{align}
Moreover, if $ \mathbf 1 _{[-1,1]}\le \widehat \psi \le \mathbf 1 _{[-2,2]} $ 
we have the following more particular estimate.  
For 
all intervals $ \omega,I$  satisfying   $\rho :=\abs{ I}\abs{ \omega } ^{-1}>1$,  and $ t>0$
\begin{equation}\label{e.zxc2}
\sum _{\substack{s\in \mathcal T(\{\omega \})\\ 
I_s\subset I}} \abs{ \ip f,  \psi _s .  } ^{2} \lesssim
(t \rho ) ^{-5}\norm \psi  ^{3 \omega }\ast f.2.^2
\qquad \text{ $ f$ supported on $ [ t  I] ^{c} $.}
\end{equation}
\end{proposition}

In the second inequality observe that we assume $\abs{ I}\abs{ \omega } ^{-1} >1$, 
so that the rectangle $ I\times \omega $ is too big to be a tile. It is important 
that on the right hand side we have both a condition on the spatial support of
 of $ f$, and in  the norm we are making a convolution with a smooth analog 
 of a Fourier projection.

\begin{proof}
The hypothesis (\ref{e.zxc}) is too strong; we are not interested in the minimal hypotheses 
here, but it is useful for this proof to observe that we only need 
\begin{equation}\label{e.zxC}
\Abs{\operatorname  \psi (x) } \lesssim (1+\abs x) ^{-5}\,,
\end{equation}
to conclude the first inequality (\ref{e.zxc0}). 

The inequality (\ref{e.zxc0}) can be seen as the assertion of the boundedness 
of the map $ f\longrightarrow \{ \ip f,   \psi_s . \mid s\in \mathcal T(\{\omega \})\}$ 
from $ L^2$ to $ \ell^2(\mathcal T (\{\omega \}))$.   It is equivalent to show that the formal dual of this 
operator is bounded, and this inequality is 
\begin{equation}\label{e.zxcdual}
\NOrm \sum _{ s\in \mathcal T(\{\omega \})} a_s \psi_s   .2. \lesssim \norm a_s.\ell^2
(\mathcal T(\{\omega \}).\,.
\end{equation}

Observe that 
\begin{align*}
\ip \psi_s ,\psi_ {s'}. \lesssim \Delta (s,s'):=
(1+\abs{ I_s} ^{-1}\abs{ c(I_s)-c(R _{s'}}) ^{-5}
\end{align*}
Estimate 
\begin{equation} \label{e.schur}
\begin{split}
\NOrm \sum _{s\in \mathbb Z} a_s \psi_s   .2.
&\le \sum _{s}\abs{ a_s} \sum _{s'} \abs{ a _{s'}} 
\Delta(s,s')
\\
&\le \norm a_s. \ell^2. \Biggl[ \sum _{s} \ABs{ \sum_{s'} \abs{ a _{s'}} 
\Delta(s,s')
} ^{2} \Biggr] ^{1/2}
\\
& \lesssim   
\norm a_s. \ell^2. \Bigl[ \sum _{s}  \sum_{s'} \abs{ a _{s'}}^2 
\Delta(s,s')
 \Bigr] ^{1/2}
\\ & \le \norm a_s .\ell^2.^2\,.
\end{split}\end{equation}
Here, we use Cauchy Schwartz, and the fact that the $ L^2$ norm dominates the $ L^1$ 
norm on probability spaces. 

\medskip
Turning to the proof of more particular assertation (\ref{e.zxc2}), 
we first note a related inequality.  Assume that 
$ \psi $ satisfies (\ref{e.zxc}). 
\begin{equation}\label{e.zxc22}
\sum _{\substack{s\in \mathcal T(\{\omega \})\\ 
I_s\subset I}} \abs{ \ip f,  \psi _s .  } ^{2} \lesssim
(t \rho ) ^{-5}\norm  f.2.^2
\qquad \text{ $ f$ supported on $ (tI) ^{c} $.}
\end{equation}
As in the statement of the Lemma, $ \rho = \abs{ I} \abs{ \omega } ^{-1}>1$. 
Here, we do not assume that $ \psi $ has compact frequency support, just 
that it has rapid spatial decay.  On the right hand side, we do not 
impose the convolution with $ \psi  ^{3 \omega }$.

For an interval $ I$ of length at 
least one, and $ t>1$, write $ \psi =\psi _0+\psi _{\infty }$ where 
$ \psi _{\infty }(x)$ is supported on $ \abs{ x}\ge \tfrac14 t \rho  $, 
equals $ \psi (x)$ on $ \abs{ x }\ge \tfrac12t \rho  $, and satisfies the estimate 
\begin{equation*}
\abs{ \psi _{\infty }(x)} \lesssim (t \rho ) ^{-10} (1+\abs{ x}) ^{-5}\,.
\end{equation*}

That is, $ \psi _{\infty } $ satisfies the inequality (\ref{e.zxC}) with constants 
that are smaller by an order of $ (t \rho ) ^{-10}$.   

Note that if $ f$ is supported on the complement of $ tI$, we have 
$ \ip f, \psi_s .=\ip f, \psi _{\infty } .$ 
for $ \lambda _s\in I$.  Thus,  (\ref{e.zxc22}) follows. 

\smallskip 

We now prove (\ref{e.zxc}) as stated.  We now assume that $ \psi $ is a Schwartz 
function satisfying $ \mathbf 1 _{[-1,1]}\le \widehat \psi \le \mathbf 1 _{[-2,2]}$. 
Then certainly, it satisfies (\ref{e.zxc}), so that (\ref{e.zxc22}) holds.  We also 
have that  for all tiles $ s\in \mathcal T(\{\omega \})$,
\begin{equation*}
 \ip f,\psi _s.=\ip \psi ^{3 \omega }\ast f,\psi .
 =\ip \psi ^{3 \omega }\ast \psi ^{3 \omega }\ast f,\psi .\,.
\end{equation*}
  Write $ \psi ^{3 \omega }\ast  \psi ^{3 \omega }\ast f=F_0+F_\infty $, where 
  $ F_0=[\psi ^{3 \omega }\ast f ]\mathbf 1 _{
\frac t 2 I}$.  
  
Then, since $ \psi $ is decreasing rapidly, we will have $ \norm F_0.2. \lesssim (t 
\rho ) ^{-10} \norm \psi ^{3 \omega }\ast  f.2.$.  Therefore, by the $ L^2$ inequality (\ref{e.zxc0}) 
\begin{equation*}
\sum _{s\in \mathcal T(\{\omega \})} \abs{ \ip F_0,  \psi _s .  } ^{2} \lesssim
(t \rho ) ^{-10}\norm  \psi ^{3 \omega }\ast  f.2.^2\,.
\end{equation*}
On the other hand, the inequality (\ref{e.zxc22}) applies to $ F _{\infty }$, so that 
\begin{equation*}
\sum _{\substack{s\in \mathcal T(\{\omega \})\\ 
I_s\subset I}} \abs{ \ip F _{\infty },  \psi _s .  } ^{2} \lesssim
(t \rho ) ^{-5}\norm  F_\infty .2.^2\,.
\end{equation*}
But certainly $ \norm F _{\infty }.2.\le \norm \psi ^{3 \omega }\ast f.2. \lesssim \norm f.2.$.
So our proof of the more particular assertation (\ref{e.zxc2}) is finished.
\end{proof}

Let us now argue that the tile operator $ \operatorname T ^{\Omega }$ maps $ L^2$ 
into itself, under the assumption that $ \Omega $ is well distributed.   
For $ \omega \in \Omega $, let $ \mathcal T(\omega )$ be the tiles in $ s\in \mathcal T(\Omega )$
 with $ \omega _s=\omega $.  
 It follows from Proposition~\ref{p.L2}, that we have the estimate 
 \begin{equation*}
\sum _{s\in \mathcal T(\omega )} \abs{ \ip f,\varphi _s.}^2 \lesssim \norm f.2.^2\,.
\end{equation*}
For a tile $ s$, we have $ \ip f ,\varphi _{s}.=
\ip \operatorname S _{2\omega }f,\varphi _s.$,
where we impose the Fourier projection onto the interval $ 2 \omega _s$ in the second 
inner product.  Thus, on the right hand side above, we can replace $ \norm f.2. ^2$ by 
$ \norm \operatorname S _{2 \omega } f.2.^2$.

Finally, the well distributed assumption implies that 
\begin{equation*}
\sum _{\omega \in \Omega } \norm \operatorname S _{2 \omega } f.2.^2 \lesssim \norm f.2.^2\,.
\end{equation*}
The boundedness of the tile operator on $ L^2$ follows.

\medskip

To prove the remaining inequalities,  we seek an appropriate endpoint estimate. 
That of $BMO$ is very useful.  Namely for $f\in L^\infty$, we show that 
 \begin{equation}\label{e.bmoT}
\norm \operatorname (T^\Omega f) ^{2}.BMO.\lesssim{}\norm f.\infty.^2\,. 
 \end{equation} 
Here, by $BMO$ we mean dyadic $BMO$, which has this definition.
 \begin{equation}\label{e.bmo} 
\norm g.BMO.=\sup_{I\in\mathbf D} \dashint_I \Abs{g-\dashint_I g}\; dx\,.
 \end{equation}
The usual definition of $BMO$ is  formed by taking a supremum over all intervals, 
not just the dyadic ones. It is a 
useful simplification for us to restrict the supremum to dyadic intervals. 
The $L^p$ inequalities for $\operatorname T^\Omega$ are deduced by an interpolation argument, which we will summarize below. 

There is a closely related notion, one that in the one parameter setting 
coincides with the $BMO$ norm.  We distinguish it here, as it is a useful distinction 
for us in the higher parameter case. 
For a map $\alpha\mid \mathbf D\longrightarrow\mathbb R$, set 
 \begin{equation} \label{e.CM} 
\norm \alpha .CM.=\sup_{J\in\mathbf D}  \abs{J}^{-1} \sum_{I\subset J}\abs{\alpha(I)}  .
 \end{equation}
``CM'' is for Carleson measure. The inequality (\ref{e.bmoT}) is, in this notation 
 \begin{equation} \label{e.CMT}
\NOrm \Bigl\{ \sum_{\substack{ s\in\mathcal T(\Omega)\\ I_s=J}} \abs{\ip f,\varphi_s.}^2\mid J\in\mathbf D\Bigr\} .CM.\lesssim{} \lVert f\rVert_\infty^2.
 \end{equation}
Or, equivalently, that we have the inequality 
\begin{equation*}
\sum _{\substack{s\in \mathcal T(\Omega )\\ I_s\subsetneq J }} \abs{\ip f,\varphi_s.}^2
\lesssim \abs{ J} \norm f.\infty .^2\,.
\end{equation*}
Notice that we can restrict the sum above to tiles $ s$ with $ I_s\subsetneq J$ 
as in the definition of $ BMO$ we are subtracting off the mean.

\begin{proof}[Proof of (\ref{e.bmoT}).]
 
Our proof  follows a familiar pattern of argument. 
Fix a function $f$ of $L^\infty$ norm one. 
We fix a dyadic interval $J$ on which we check the $BMO$ norm.  
We write $ f=\sum _{k=1} ^{\infty } g_k$, where 
$g_1=f\ind {2J}$, and  
 \begin{equation*}
g_k=f\ind {2^{k}J-2^{k-1} J },\qquad k>1.
 \end{equation*} 
 The bound below follows from the $ L^2 $ bound on the tile operator. 
 \begin{equation*}
 \sigma (k):=\sum_{\substack{s\in \mathcal T(\Omega )\\ I_s\subsetneq J }}
\abs{\ip g_k,\varphi_s.}^2\lesssim{}\lVert g_k\rVert_2^2\lesssim{} 2 ^{k} \abs J
 \end{equation*}
  For $ k>5$, we will use the more particular estimate (\ref{e.zxc2}) 
  to verify that 
\begin{align} \label{e.sigmak}
\sigma (k)^2 :=  \sum_{\substack{s\in \mathcal T(\Omega )\\ I_s\subsetneq J }}
\abs{\ip g_k,\varphi_s.}^2\lesssim{} 2 ^{-4k}\lVert g_k\rVert_2^2
\lesssim{} 2 ^{-4k}  \abs{ J}\,.
\end{align} 
Yet, to apply  (\ref{e.zxc2}) we need to restrict attention to 
a single frequency interval $ \omega $, which we do here.  
\begin{equation*}
\sum_{\substack{s\in \mathcal T(\Omega )\\ I_s\subsetneq J\,, \omega _s=\omega  }}
\abs{\ip g_k,\varphi_s.}^2\lesssim{} 2 ^{-10k} \norm \varphi ^{3 \omega }*g_k.2.^2, 
\qquad \omega \in \Omega \,.
\end{equation*}
This is summed over $ \omega \in \Omega $, using the estimate
\begin{equation*}
\sum _{\omega \in \Omega } \norm \varphi ^{3 \omega }*g_k.2.^2 
\lesssim \norm g.2.^2 \lesssim 2 ^k \abs{ J}
\end{equation*}
to prove (\ref{e.sigmak}).

The inequality (\ref{e.sigmak}) is summed over $ k$
in the following way to  finish the proof of the $ BMO $ estimate, 
(\ref{e.bmoT}). 
\begin{equation}\label{e.summed}
\begin{split}
\sum_{\substack{s\in \mathcal T(\Omega )\\ I_s\subsetneq J }}
\abs{\ip f,\varphi_s.}^2 &
=\sum_{\substack{s\in \mathcal T(\Omega )\\ I_s\subsetneq J }}
\ABs{ \sum _{k=1} ^{\infty } k ^{-1} \cdot k ^{1} \cdot \ip g_k,\varphi _s.} ^{2}
\\
&\lesssim \sum _{k=1} ^\infty k^2 \sigma (k)^2 \lesssim \abs{ J}\,.
\end{split}
\end{equation}

\end{proof}

We discuss how to derive the $ L^p$ inequalities from the $ L^2$ estimate and 
the $ L^\infty \longrightarrow \textup{BMO}$ estimate.

The method used by Rubio de Francia \cite {rubio}, to use our notation, was to prove 
the inequality  $ [(\operatorname T ^{\Omega }f)^2 ]^\sharp \lesssim \operatorname M \abs{ f}^2$, 
where $ g ^{\sharp}$ is the (dyadic) sharp function defined by 
\begin{equation*}
g ^{\sharp}(x)=\sup _{\substack{x\in I\\ I \in \mathbf D}} \dashint _{I} 
\Abs{ g(y)-\dashint _{I} g(z)\; dz}\; dy 
\end{equation*}
One has the inequality $ \norm g ^{\sharp}.p. \lesssim \norm g.p.$ for $ 1<p<\infty $.  
The proof we have given can be reorganized to prove this estimate.

We have not presented this argument since the sharp function does not permit 
a good extension to the case of higher parameters, which we discuss in the next section. 
On the other hand, a proof of the (standard) interpolation result between $ L^p$ and $ BMO$ 
\cite {MR928802} is 
based upon the John Nirenberg inequality, Lemma~\ref{l.jn} below;  a proof based upon  
this inequality does extend to higher parameters.  We present this argument now. 

One formulation of the inequality of F.~John and L.~Nirenberg is 

 \begin{lemma}\label{l.jn}  
For each $1<p<\infty$, we have the estimate below valid for all dyadic intervals $J$, 
 \begin{equation*}
\NOrm \sum_{I\subset J} \frac{\alpha(I)}{\abs I}\ind I .p.\lesssim{} \norm \alpha .CM.\abs{J}^{1/p}.
 \end{equation*}
The implied constant depends only on $ p$.
 \end{lemma}

\begin{proof} It suffices to prove the inequality for $p$ an integer, as the remaining values of $p$ are available 
by H\"older's inequality.  The case of $p=1$ is the definition of the Carleson measure norm.  Assuming the inequality for 
$p$, consider 
 \begin{align*} 
\int_J \Bigl[ \sum_{I\subset J} \frac{\alpha(I)}{\abs I}\ind I \Bigr]^{p+1} \; dx{}\le{}&2 
	\sum_{J'\subset J} \frac{\abs{\alpha(J')}}{\abs{J'}}\int_{J'} \Bigl[ \sum_{I\subset J'} \frac{\alpha(I)}{\abs I}\ind I \Bigr]^{p} \; dx
\\{}\lesssim{}& \lVert \alpha\rVert_{CM}^p \sum_{J'\subset J} \abs{\alpha(J')} 
\\{}\lesssim{}& \lVert \alpha\rVert_{CM}^{p+1}\abs{J}.
 \end{align*}
Notice that we are strongly using the grid property of the dyadic intervals, namely that for $I,J\in\mathbf D$ we have 
$I\cap J\in\{\emptyset,I,J\}$. 

For an alternate proof, see Lemma~\ref{l.jnproduct} below.
\end{proof}

We  prove the following operator for the tile operator $ \operatorname T ^{\Omega }$:
 \begin{equation} \label{e.restricted}
\norm \operatorname T^\Omega\ind F .p.\lesssim{} \abs{F}^{1/p}, \qquad 2<p<\infty,
 \end{equation} 
for all sets $F\subset \mathbb R$ of finite measure. This is the restricted 
strong type inequality on $ L^p$ for the tile operator---that is we only prove the $ L^p$ 
estimate for indicator functions.

The $L^p$ inequality above is obtained by considering 
subsets of tiles, $\mathcal T\subset\mathcal T(\Omega)$, for which we will need  the notation 
 \begin{equation*}
\operatorname T^{\mathcal T}\ind F:=\Bigl[ \sum_{s\in\mathcal T} \frac{\abs{\ip \ind F,\varphi_s.}^2 }{\abs{I_s}} \ind {I_s} \Bigr]^{1/2}
 \end{equation*}
As well, take $\sh {\mathcal T}:=\bigcup_{s\in\mathcal T}I_s$ to be the {\em shadow of $\mathcal T$.}

 The critical step is to decompose $\mathcal T(\Omega)$ into 
 subsets $\mathcal T_k$ for which 
 \begin{equation} \label {e.tk}   
\norm (\operatorname T^{\mathcal T_k}\ind F )^2 .BMO.\lesssim{} 2^{-2k},\qquad  \abs{\sh {\mathcal T_k}}\lesssim{}2^{2k}\abs F, \qquad k\ge1.
 \end{equation}
We have already seen that  the $BMO$ norm is bounded, so we need only consider $k\ge1$ above. 
Then, by the John--Nirenberg inequality, 
 \begin{equation*}
\norm  \operatorname T^{\mathcal T_k}\ind F .p.\lesssim 2^{-k (1-2/p)}\abs F ^{1/p}.
 \end{equation*}
This is summable in $k$ for $p>2$.  

The decomposition (\ref{e.tk}) follows from this claim.  Suppose that $\mathcal T\subset \mathcal T(\Omega)$ satisfies 
 \begin{equation*} 
 \norm (\operatorname T^{\mathcal T}\ind F )^2.BMO.\lesssim\beta
  \end{equation*}
 We show how to write it as a union of $\mathcal T_{\text{big}}$ and $\mathcal T_{\text{small}}$ where 
  \begin{equation*}
\norm ( \operatorname T^{\mathcal T_{\text{small}}}\ind F)^2 .BMO.\lesssim\tfrac\beta4
\,,\qquad 
\abs{ \sh{\mathcal T_{\text{big}}} }\lesssim{}\beta^{-1} \abs{ F}\,.
  \end{equation*}

 The decomposition is achieved in a recursive fashion.  Initialize 
  \begin{equation*}
 \mathbf J:=\emptyset \quad 
 \mathcal T_{\text{big}}:=\emptyset,\quad \mathcal T_{\text{small}}:=\emptyset, \quad
 \mathcal T_{\text{stock}}:=\mathcal T.
  \end{equation*}
 While  $\norm (\operatorname T^{\mathcal T_{\text{stock}}}\ind F )^2.BMO.\ge\tfrac\beta4$,
 there is a maximal dyadic interval $J\in\mathbf D$ for which 
  \begin{equation*}
 \sum_{\substack{s\in\mathcal T_{\text{stock}} \\ I_s\subset J} }
 \abs{\ip \ind F,\varphi_s.}^2\ge\tfrac{\beta}4 \abs J.
  \end{equation*}
 Update 
 \begin{gather*}
\mathbf J:=\mathbf J\cup\{J\},\quad 
\mathcal T_{\text{big}}:=\mathcal T_{\text{big}}\cup \{s\in \mathcal T_{\text{stock}}\mid I_s\subset J\},
\\
\mathcal T_{\text{stock}}:=\mathcal T_{\text{stock}}-\{s\in \mathcal T_{\text{stock}}\mid I_s\subset J\} . 
 \end{gather*}
 Upon completion of the While loop, update $\mathcal T_{\text{small}}:=\mathcal T_{\text{stock}}$ 
 and return the values of $\mathcal T_{\text{big}}$ and $\mathcal T_{\text{small}}$. 
 
 Observe that by the $ L^2$ bound for the tile operator we have 
  \begin{align*}
\beta \abs{\sh {\mathcal T_{\text{big}}}}\lesssim{}& \beta \sum_{J\in\mathbf J} \abs J 
\\{}\lesssim{}& \sum_{s\in \mathcal T_{\text{big}}} \abs{\ip \mathbf 1 _{F},\varphi_s.}^2 
\\{}\lesssim{}& \abs F .
 \end{align*}
This completes the proof of (\ref{e.tk}).  Our discussion of the restricted strong 
type inequality is complete.

\section{The Case of  Higher Dimensions}

We give the proof of \thm.lp/ in higher dimensions.  
The tensor product structure permits us to adapt many of
the arguments of the one dimensional case.  (Some arguments are far less trivial to adapt 
however.) For 
instance, one can apply the classical Littlewood Paley inequality in each variable separately.  This would yield 
a particular instance of a Littlewood Paley inequality in higher dimensions.  Namely, for all dimensions $d$, 
 \begin{equation} \label{e.higherlp}
\norm \operatorname S^{\Delta^d}f.p.\simeq \norm f.p.,\qquad 1<p<\infty, 
 \end{equation}
where $\Delta^d=\bigotimes_1^d \Delta$ is the $d$--fold tensor product of the lacunary intervals $\Delta$, as in \thm.lp/.  
  
Considerations of this type apply to many of the arguments made in the one dimensional case of \thm.lp/.  
In particular the definition of well distributed, and the \l.wellsuffice/ continues to hold in the higher dimensional 
setting.    

As before, the well distributed assumption permits one to define a
``smooth'' square function that is clearly bounded on $L^2$. 
We again choose to replace a convolution square function with an appropriate 
tile operator.

The definition of  the smooth square function---and of  tiles---requires a little more care.  For positive quantities 
$t=(t_1,\ldots,t_d)$, dilation operators  are given by 
 \begin{equation*}
\Dil p t f(x_1,\ldots,x_d)= \Bigl[\prod_{j=1}^d t_j^{-1/p}\Bigr] 
f(x_1/t_1,\ldots,x_j/t_d),\qquad 0\leq{}p\leq\infty\,,
 \end{equation*}
with the normalization chosen to preserve the $ L^p$ norm of $ f$. 

A rectangle is a product of intervals in the standard basis.  Writing a rectangle 
as $R=R_{(1)}\times \cdots \times R_{(d)}$,  we extend the definition of the dilation operators 
in the following way. 
 \begin{equation*}
\Dil p R := \Tr {c(R)}\Dil p {(\abs {R_{(1)}},\ldots,\abs{R_{(d)}})} 
 \end{equation*}

 For a Schwartz function $\varphi$ on $\mathbb R^d$, satisfying 
  \begin{equation*}
 \ind {[-1/2,1/2]^d}\le\widehat\varphi\le{}\ind {[-1,1]^d}
  \end{equation*}
 we set 
  \begin{equation} \label{e.zvfzw}
 \varphi^\omega=\Mod {c(\omega)} \Dil 1 { (\abs{\omega_{(1)}}^{-1},\ldots,\abs{\omega_{(d)}}^{-1})} \varphi 
  \end{equation}
  For a collection of well distributed rectangles $\Omega$, we should show that the  inequality (\ref{e.G}) holds.

\begin{figure}
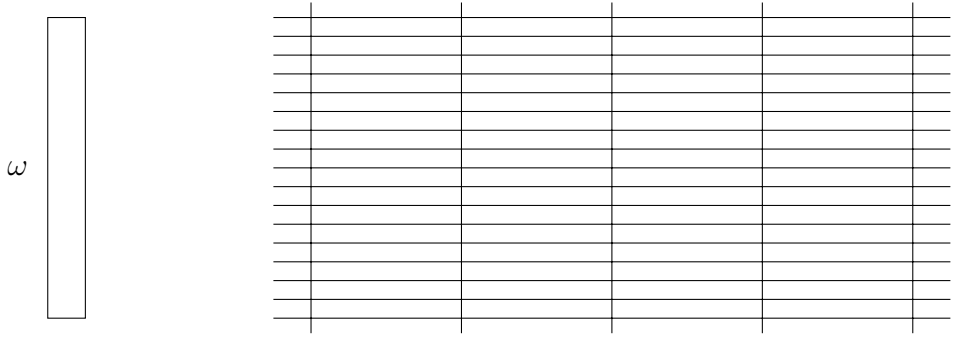
  \label{f.2dTiles}
\begin{center}
\begin{pgfpicture}{0cm}{0cm}{8cm}{4cm}
\begin{pgftranslate}{\pgfxy(4,2)}  
\pgfgrid[stepx=2cm,stepy=.25cm]{\pgfxy(-2.5,-2.2)}{\pgfxy(6.5,2.2)}
\end{pgftranslate}
\pgfrect[stroke]{\pgfxy(-1.5,0)}{\pgfxy(.5,4)}
\pgfputat{\pgfxy(-1.9,2)}{\pgfbox[center,center]{$\omega$}}  
\end{pgfpicture}
\end{center}
\caption{Dual rectangles in the plane: An example frequency rectangle $\omega$ on the left, and dual dyadic 
rectangles on right. }
\end{figure}

  We substitute the smooth convolution square function for a sum over tiles.   
  Say that $R\times\omega$ is a {\em  tile} if  both $\omega$ and $R$ are rectangles and 
  for all $1\le j\le d$, 
  $1\le\abs {\omega_{(j)}}\cdot\abs{R_{(j)}}<2$, and $R_{(j)}$ is a dyadic interval. 
  Thus, we are requiring that $ \omega $ and $ R$ be dual in each coordinate seperately. 
  In this instance, we  refer to $ \omega $ and $ R$ as dual rectangles.
    See Figure~\ref{f.2dTiles}.

 Write $s=R_s\times \omega_s$. 
  As before, let $\mathcal T(\Omega)$ be the set of all tiles $s$ such that $\omega_s\in\Omega$. 
Define functions adapted to tiles and a tile operator by 
   \begin{gather*} 
  \varphi_s={} \Mod {c(\omega_s)} \Dil 2 {R_s} \varphi 
   \\
 \operatorname T^\Omega f=\Bigl[\sum_{s\in\mathcal T(\Omega)} \frac{ \abs{\ip f,\varphi_s.}^2}{{\abs{R_s}}} \ind {R_s} \Bigr]^{1/2}
  \end{gather*}
 The main point is to establish the boundedness of this operator on $ L^2$ and 
 an appropriate endpoint estimate.  The analog of Proposition~\ref{p.L2} is 
 in this setting

\begin{proposition}\label{p.Td}   Assume only that the function 
$ \varphi $ satisfies 
\begin{equation}\label{e.Tdphi}
\abs{ \varphi (x)} \lesssim (1+\abs{ x}) ^{-20d}\,.
\end{equation}
Let $ \Omega =\{\omega \}$.  Then, we have the estimate 
\begin{equation}\label{e.TdT}
\sum _{s\in \mathcal T(\{\omega\})  } \abs{ \ip f,\varphi _s.} ^{2} \lesssim \norm f.2.^2\,.
\end{equation}
Now assume that $ \varphi $ is a smooth Schwartz function satisfying 
$ \mathbf 1 _{[-1,1]^d}\le \widehat  \varphi \le \mathbf 1 _{[-2,2]^d} $. 
For a subset $U \subset \mathbb R ^{d}$ of finite measure, $ 0<a<1$, and  
a function  $ f$ supported on the complement of 
$ \{ \operatorname M \mathbf 1 _{U} > a \}$, we have the estimate 
\begin{equation}\label{e.TdU}
  \sum_{\substack{R_s\subset U \\ \omega_s=\omega}}\abs{\ip f,\varphi_s. }^2
  \lesssim{}a^{15d} \norm  \varphi ^{3 \omega }\ast f.2.^2
\end{equation}
\end{proposition}

The more particular assertation (\ref{e.TdU})  has a far more complicated form than 
in the one dimensional setting.  That is because when we turn to the endpoint estimate, 
it is a Carleson measure condition; this condition is far more subtle, in that 
it requires  testing the measure  against arbitrary sets, instead of just intervals, or rectangles.

\begin{proof}

The hypothesis (\ref{e.Tdphi}) is more than enough to conclude (\ref{e.TdT}). 
We need only assume 
\begin{equation}\label{e.TdphI}
\abs{ \varphi (x)} \lesssim (1+\abs{ x}) ^{-5d}\,.
\end{equation}

After taking an appropriate dilation and modulation,
we can assume that $ \omega =[-\tfrac12,\tfrac12] ^{d}$.   We view the inequality 
(\ref{e.TdT}) as the boundedness of the linear map $ f \longrightarrow \{ \ip 
f,\varphi _s. \mid s\in \mathcal T(\{\omega \})\} $ from $ L^2(\mathbb R ^{d} ) $ 
into $ \ell ^2 (\mathbb Z ^{d})$.  We then prove that the dual to this operator 
is bounded, that is we verify the inequality 
\begin{equation*}
\NOrm \sum _{s\in \mathcal T(\{\omega \})} a _{s} \varphi _s .2. 
\lesssim \norm a_s . \ell ^{2} (\mathcal T(\{\omega \}).
\end{equation*}

Observe that 
\begin{equation*}
\abs{ \ip \varphi _s, \varphi _{s'}.} \lesssim (1+\text{dist}(R_s, R_ {s'})) ^{-5d}\,,
\qquad s,s'\in \mathcal T(\{\omega \})\,.
\end{equation*}
The remaining steps of the proof are a modification of (\ref{e.schur}). 

\medskip 
As in the one dimensional setting, the more particular assertation is proved in 
two stages.  First we assume only that the function $ \varphi $ satisfy 
(\ref{e.Tdphi}), and prove 
\begin{equation}\label{e.TdUX}
  \sum_{\substack{R_s\subset U \\ \omega_s=\omega}}\abs{\ip f,\varphi_s. }^2
  \lesssim{}a^{15d} \norm  f.2.^2
\end{equation}
for functions $ f$ supported on the complement of 
$ \{ \operatorname M \mathbf 1 _{U} > a \}$.  

 Take $ \widetilde \varphi (x)$ to be a function 
which equals $ a ^{-15d} \varphi (x) $ provided $ \abs{ x}\ge \tfrac 2 a$. 
With this,  $ \varphi $ satisfies (\ref{e.TdphI}) with a constant independent of $ a$. 

For a subset $ U\subset \mathbb R ^{d}$ of finite measure, and function 
$ f$ supported on the complement of $ \{ \operatorname M \mathbf 1 _{U} > a^d\}$, 
and tile $ s$ with $ R_s\subset U$, we have $ a ^{-15d}\ip f ,\varphi _s.=\ip f,\widetilde \varphi_s. $.
Thus, (\ref{e.TdUX}) follows from the $ L^2$ estimate we have already proved. 

\smallskip 

We can then prove the assertation of the Lemma. Take $ \varphi $, $ f \in L^2$, 
$ 0<a<1$, and $ U\subset \mathbb R ^{d}$ as in (\ref{e.TdU}).  
Then, for all tiles $ s=R_s\times \omega $, we have 
\begin{equation*}
 \ip f,\varphi _s.=
\ip \varphi ^{3 \omega } \ast f,\varphi _s.
\end{equation*}
We write $\varphi ^{3 \omega }\ast f =F_0+F _{\infty }$, where 
\begin{equation*}
F_0=[  \varphi ^{3 \omega }\ast f ]\mathbf 1 _{ \{\operatorname M \mathbf 1 _{U} > 2a\}}\,.
\end{equation*}

The rapid decay of $ \varphi $, with the fact 
about the support of $ f$, show that 
$ \norm F_0.2. \lesssim a ^{15d} \norm \varphi ^{3 \omega } \ast  f.2.$.  
Thus, the estimate below follows from the $ L^2$ inequality (\ref{e.Tdphi}):
\begin{equation*}
\sum_{\substack{R_s\subset U \\ \omega_s=\omega}}\abs{\ip F_0,\varphi_s. }^2
  \lesssim{}a^{15d} \norm \varphi ^{3 \omega } \ast  f.2.^2\,.
\end{equation*}
As for the term $ F_\infty  $, we use the  estimate 
(\ref{e.TdUX}) to see that 
\begin{equation*}
\sum_{\substack{R_s\subset U \\ \omega_s=\omega}}\abs{\ip F_\infty ,\varphi_s. }^2
  \lesssim{}a^{15d} \norm \varphi ^{3 \omega } \ast  f.2.^2
\end{equation*}
This completes our proof of (\ref{e.TdU}).

\end{proof}

We can now prove the $ L^2 $ boundedness of the square function.  Using the 
well distributed assumption and (\ref{e.TdT}), we can estimate 
\begin{align*}
\sum _{s\in \mathcal T(\Omega )} \abs{ \ip f,\varphi _s.} ^2 
& =\sum _{\omega \in \Omega } 
\sum _{s\in \mathcal T(\{\omega \})} \abs{ \ip  S _{2\omega }f,\varphi _s.} ^2
\\ & \lesssim  \sum _{\omega \in \Omega } \norm S _{2\omega }f  .2.^2
\\ & \lesssim \norm f.2.^2\,. 
\end{align*}

The endpoint estimate we seek is phrased this way.  For all subsets 
$ U\subset \mathbb R ^{d}$ of finite measure, and functions $ f$ of $ L^ \infty $ norm one, 
\begin{equation}\label{e.BMOT}
\abs{ U} ^{-1} \sum _{\substack{s\in \mathcal T(\Omega )\\ R_s \subset U }}
\abs{ \ip f, \varphi _s.} ^{2} \lesssim 1. 
\end{equation}
Using the notation of (\ref{e.productCM}), this inequality is equivalent to 
\begin{equation*}
\NOrm \Bigl\{  \sum _{\substack{s\in \mathcal T\\  R_s=R} } 
\frac{ \abs{ \ip f, \varphi _s.} ^{2} } {\abs{ R}} \mathbf 1 _{R} \mid R  \in 
\mathbf D ^{d} \Bigr\} . CM. \lesssim \norm f.\infty .\,.
\end{equation*}

 Write  $f=\sum _{k=1} ^{\infty } g_k $ where 
   \begin{gather*} 
  g_1=f\ind { \{ \operatorname  M \ind U \ge \tfrac1{2} \} } , 
  \\
  g_k=f\ind { \{ 2^{-k}\le \operatorname   M \ind U \le 2^{-(k-1)}\}}, \qquad k>1.
   \end{gather*}
   Using the boundedness of the maximal function on e.g.~$ L^2$, and the $ L^2$ 
   boundedness of the tile operator, we have 
   \begin{equation*}
\abs{ U} ^{-1} \sum _{\substack{s\in \mathcal T(\Omega )\\ R_s \subset U }}
\abs{ \ip g_k, \varphi _s.} ^{2} \lesssim  \abs{ U} ^{-1}\norm f_k.2.^2 \lesssim 2 ^{2k}\,. 
\end{equation*}

For the terms arising from $ g_k$, with $ k\ge5$,  we can use (\ref{e.TdU}) with 
$ a=2 ^{-k/d}$ to see that 
\begin{align*}
\sum _{\substack{s\in \mathcal T(\Omega )\\ R_s \subset U }}
\abs{ \ip g_k, \varphi _s.} ^{2} 
&=\sum _{\omega \in \Omega } \sum _{\substack{s\in \mathcal T(\Omega )\\ \omega = \omega _s\,,\ 
R_s \subset U }}
\abs{ \ip g_k, \varphi _s.} ^{2} 
\\
& \lesssim 2 ^{-10d}  \sum _{\omega \in \Omega } 
\norm \varphi ^{3 \omega } \ast g_k.2.^2 
\\
& \lesssim 2 ^{-10k} \norm g_k.2.^2 
\\ &\lesssim 2 ^{-8k} \abs{ U}. 
\end{align*}
Here, we have used the fact that the strong maximal function is bounded on $ L^2$. 
The conclusion of the proof of (\ref{e.BMOT}) then follows the lines of (\ref{e.summed}).

\medskip 

To deduce the $ L^p$ inequalities, one can again appeal to interpolation.  Alternatively, 
the restricted strong type inequality can be proved directly using the 
John Nirenberg inequality for the product Carleson measure.  This inequality is 
recalled in the next section, and argument is formally quite simliar to the one we 
gave for one dimension.  Details are omitted.

\subsubsection*{Carleson Measures in the Product Setting}

The subject of Carleson measures are central to the subject of product $ BMO$, as 
discovered by S.-Y.~Chang and R.~Fefferman \cites{MR86g:42038
,
MR82a:32009}.

A definition can be phrased in terms of  maps $ \alpha $ from the dyadic rectangles $\mathbf D^d$ of $\mathbb R^d$.
 This norm is 
  \begin{equation} \label{e.productCM}
\norm \alpha .CM.={}\sup_{U\subset \mathbb R^d} \abs{U}^{-1}\sum_{R\subset U} \alpha(R).
 \end{equation}
What is most important is that the supremum is taken over all sets $U\subset\mathbb R^d$ of finite measure.  It would of course 
be most natural to restrict the supremum to rectangles, and while this is not an adequate definition, it nevertheless plays an 
important role in the theory.   See the Lemma of Journ\'e \cite{MR87g:42028}, as well as the 
survey of Journel's Lemma of Cabrelli, Lacey, Molter, and Pipher \cite{journesurvey}.

Of importance here is the analog of the John--Nirenberg inequality in this setting. 

 \begin{lemma}\label{l.jnproduct}  We have the inequality below, valid for all sets $U$ of finite measure. 
 \begin{equation*}
\NOrm \sum_{R\subset U} \frac {\alpha(R)}{\abs R}\ind R .p.\lesssim{} \norm \alpha .CM. \abs{U}^{1/p},\qquad 1<p<\infty.
 \end{equation*}
 \end{lemma}

\begin{proof} 
We use the duality argument of Chang and Fefferman \cite{MR82a:32009}.
Let $\norm \alpha.CM.=1$. 
Define 
 \begin{equation*}
F_V:={}
\sum_{R\subset V} \frac {\alpha(R)}{\abs R}\ind R 
 \end{equation*}
We shall show that for all $U$, there is a set $V$ satisfying $\abs{V}<\frac12\abs U$ for which 
 \begin{equation} \label{e.UV}
\norm F_U .p.\lesssim{} \abs{U}^{1/p}+\norm F_V .p.
 \end{equation}
Clearly, inductive application of this inequality will prove our Lemma.

The argument for (\ref{e.UV}) is by duality.  Thus, for a given $1<p<\infty$, and conjugate index $p'$, take $g\in
L^{p'}$ of norm one so that $\norm F_U.p.=\ip F_U,g.$.  Set 
 \begin{equation*}
V=\{\operatorname  Mg >K \abs {U}^{-1/p'}\}
 \end{equation*}
where $M$ is the strong maximal function and $K$ is sufficiently large so that $\abs V<\frac12\abs U$.  Then, 
 \begin{equation*}
\ip F_U,g.=\sum_{\substack{R\subset U\\R\not\subset V}} \alpha(R)\dashint_R g \; dx+\ip F_V ,g.
 \end{equation*}
The second term is at most $\norm F_V.p.$ by H\"older's inequality.  For the first term, note that the 
average of $g$ over $R$ can be at most $K\abs{U}^{-1/p'}$.  So by the definition of Carleson measure norm, it is at most 
 \begin{equation*}
\sum_{\substack{R\subset U\\R\not\subset V}} \alpha(R)\dashint_R g \; dx\lesssim{}\abs{U}^{-1/p'}\sum_{R\subset U}\alpha(R)\lesssim{}\abs{U}^{1/p},
 \end{equation*}
as required by (\ref{e.UV}). 

\end{proof}

\section{Implications for Multipliers} 

Let us consider a bounded function $m$, and define 
  \begin{equation*}
 \operatorname A_m f(x):=\int m(\xi)\widehat f (\xi)\; d\xi.
  \end{equation*}
 This is the multiplier operator given by $m$, and the Plancherel equality 
implies that  the operator  norm of $A$ on $ L^2$ is given by $\norm m.\infty.$.  It is 
 of significant interest to have a description of the the norm of $A$ as an operator on $L^p$ only in terms of properties of the 
function $m$.

 Littlewood--Paley inequalities have implications here, as is recognized through the proof of the classical Marcinciewcz Theorem.  
 Coifman, Rubio de Francia and Semmes \cite{MR89e:42009} found a beautiful extension of this classical Theorem with a  proof that  is a pleasing application of 
 Rubio de Francia's inequality.  We work first in one dimension.  
 To state it, for an interval $[a,b]$, and index $0<q<\infty$, we set the $q$ variation norm of $m$ on the interval $[a,b]$ to be 
  \begin{equation} \label{e.VAR}
 \norm m. \operatorname{Var}_q([a,b]).:=\sup\Bigl\{ \Bigl[\sum_{k=1}^{K} \abs{m(\xi_{k+1})-m(\xi_k)}^q \Bigr]^{1/q} \Bigr\} 
 \end{equation}
where the supremum is over all finite sequences $a=\xi_0<\xi_1<\xi_2<\ldots<\xi_{K+1}=b$. 
Set $\norm v.V_q([a,b]).:=\norm m .L^\infty([a,b]).+\norm m. \operatorname{Var}_q([a,b]).$
 Note that if $q=1$, this norm coincides with the 
classical bounded variation norm.

 \begin{theorem}\label{t.crs}   Suppose that $1<p,q<\infty$, satisfying $\abs{\frac 12-\frac1p}<\frac1q$.  Then for all functions $m\in L^{\infty}(\mathbb R)$, we have 
 \begin{equation*}
\norm \operatorname A_m .p.\lesssim{} \sup_{I\in\mathbf D} \norm m .V_q(I).
 \end{equation*}
 \end{theorem} 

Note that the right hand side is a supremum over the Littlewood--Paley intervals $I\in \mathbf D$. The Theorem above is 
as in the Marcinciewcz Theorem, provided one takes $q=1$.
But the Theorem of Coifman, Rubio de Francia and Semmes states that even for the much rougher case of $q=2$, the right hand side is an upper bound for all $L^p$ operator norms of the multiplier norm $\operatorname A_m$.  In addition, as $q$ increases to infinity, the $V_q$ norms approach that of $L^\infty$, which is the correct estimate for the multiplier norm at $p=2$.

 \subsection{Proof of \thm.crs/} 
 
  The first Lemma in the proof  is a transparent display of the usefulness of the Littlewood Paley inequalities in decoupling scales.  
 
  \begin{lemma}\label{l.decouple}  Suppose that the multiplier $m$ is of the form $m=\sum_{\omega\in\mathbf D}a_\omega \ind \omega $, for a sequence of reals $a_\omega$.  Then, 
  \begin{equation*}
 \norm \operatorname A_m .p.\lesssim{} \norm a_\omega .\ell^\infty(\mathbf D).,\qquad 1<p<\infty.  
  \end{equation*}
 Suppose that for an integer $n$, that $\mathbf D_n$ is a partition of $\mathbb R$ that refines the 
 partition $\mathbf D$, and partitions each $\omega\in\mathbf D$ into at most  $n$ subintervals.    Consider a multiplier of the form 
  \begin{equation*}
 m=\sum_{\omega\in\mathbf D_n}a_\omega\ind \omega .
  \end{equation*}
For $\abs{\frac1p-\frac12}<\frac1q$, we  have 
  \begin{equation} \label{e.AmDn} 
 \norm \operatorname A_m .p.\lesssim{} n^{1/q}\norm a_\omega.\ell^\infty(\mathbf D_n). 
  \end{equation} 
  \end{lemma}

 \begin{proof}  In the first claim, for each $\omega\in\mathbf D$, we have $ \operatorname S_\omega \operatorname A_m=a_\omega \operatorname S_\omega$, so that 
 for any $f\in L^p$, we have by the Littlewood Paley inequalities 
  \begin{align*}
 \norm \operatorname A_m f.p.\simeq{}& \norm S^{\mathbf D} \operatorname A_m f.p.
 \\{}\simeq{}&  \NOrm \Bigl[ \sum_{\omega\in\mathbf D}\abs{a_\omega}^2\abs{ \operatorname S_\omega f}^2\Bigr]^{1/2} .p. 
 \\{}\lesssim{}&\norm a_\omega .\ell^\infty(\mathbf D). \norm S^{\mathbf D} f.p..
  \end{align*}
  
\medskip 

 The proof of (\ref{e.AmDn}) is by interpolation. Let us presume that $\norm   a_\omega.\ell^\infty(\mathbf D_n).=1$. We certainly have  $\norm \operatorname A_m .2.=1$.  On the other hand, with an eye towards applying the classical Littlewood Paley inequality  and Rubio de Francia's extension of it, for each $\omega\in\mathbf D$, we have 
  \begin{equation*}
 \abs{\operatorname S_\omega \operatorname A_m f}\le{}n^{1/2}\Bigl[ \sum_{\substack{\omega'\in\mathbf D_n \\ \omega'\subset\omega}} \abs{ \operatorname S_{\omega'}f}^2 \Bigr]^{1/2}. 
  \end{equation*}
 Therefore, we may estimate for any $2<r<\infty$, 
  \begin{align} \label{e.Amr}
 \norm \operatorname A_m f .r.\lesssim{}& \norm \operatorname S^{\mathbf D} \operatorname A_m f.r.
 \\{}\lesssim{}& n^{1/2} \norm \operatorname  S^{\mathbf D_n} f. r. \nonumber
 \\{}\lesssim{}& n^{1/2} \norm f.r.\nonumber
  \end{align}
 To conclude (\ref{e.AmDn}), let us first note the useful principle that $\norm \operatorname A_m.p.=\norm \operatorname A_m.p'.$, where $p'$ is the conjugate index.  So we can take $p>2$.  For the choice of $\frac12-\frac1p<\frac1q$, take a value of $r$ that is very large,  in fact 
 \begin{equation*}
\frac12-\frac1p<\frac 1r<\frac1q
 \end{equation*}
 and interpolate (\ref{e.Amr}) with the $L^2$ bound. 
 \end{proof}

 \medskip 
 Since our last inequality is so close in form to the Theorem we wish to prove, the most expedient thing to do is to note a 
 slightly technical lemma about functions in the $V_q$ class. 
 
  \begin{lemma}\label{l.Vq}  If $m\in V_q(I)$ is of norm one, we can choose partitions $\Pi_j$, $j\in\mathbb N$, of $I$ into at most $2^j$ subintervals
 and functions $m_j$ that are measurable with respect to $\Pi_j$, so that 
  \begin{equation*}
 m=\sum_j m_j,\qquad \norm m_j.\infty.\lesssim{}2^{-j}. 
  \end{equation*}
 \end{lemma} 

\begin{proof}    Let $\mathbf P_j=\{(k2^{-j},(k+1)2^{-j}] \mid 0\le{}k<2^j\}$  be  the standard partition of $[0,1]$ into intervals of length $2^{-j}$.  Consider the function $\mu\mid I=[a,b]\longrightarrow[0,1]$ given by 
 \begin{equation*}
\mu(x):=\lVert m\rVert_{\operatorname{Var}_q([a,x])}^q.
 \end{equation*}
This function is monotone, non--decreasing, hence has a well defined inverse function.  
Define $\Pi_j=\mu^{-1}(\mathbf P_j)$.  
We define the functions $m_j$ so that 
 \begin{equation*}
\sum_{k=1}^j m_j=\sum_{\omega\in\Pi_j} \ind \omega \dashint_\omega m \; d\xi.
 \end{equation*}
That is, the $m_j$ are taken to be a martingale difference sequence with respect to the increasing sigma fields $\Pi_j$.  
Thus, it is clear that $m=\sum m_j$.  The bound on the $L^\infty$ norm of the $m_j$ is easy to deduce from the definitions. 
\end{proof}

We can prove the \thm.crs/ as follows.  For $\frac12-\frac1p<\frac1r<\frac 1q$, and $m$ such that 
 \begin{equation*}
\sup_{\omega\in\mathbf D} \norm m\ind \omega. V_q(\omega).\le{}1,
 \end{equation*}
we apply \l.Vq/ and (\ref{e.AmDn}) to each $m\ind \omega$ to conclude that we can write $m=\sum_j m_j$, 
so that $m_j$ is a multiplier  satisfying $\norm A_{m_j}.p.\lesssim 2^{j/r-j/q}$.  But this estimate is summable in 
$j$, and so completes the proof of the Theorem.
 
 \subsection{The Higher Dimensional Form}

 The extension of the theorem above to higher dimensions was made by Q.~Xu \cite{MR98b:42024}.  His point of view was to take 
 an inductive and vector valued approach.  Some of his ideas were motivated by prior work of 
 G.~Pisier and Q.~Xu \cites{pxu1,pxu2} in which interesting applications of $q$--variation spaces are made.   
  
  The definition of the $q$ variation in higher dimensions is done inductively.  For a function 
  $m\mid \mathbb R^d\longrightarrow\mathbb C$, define difference operators by 
  \begin{equation*}
 \operatorname{Diff}(m,k,h,x)=m(x+he_k)-m(x),\qquad 1\le k\le d,
  \end{equation*}
 where $e_k$ is the $k$th coordinate vector.  For a rectangle $R=\prod_{k=1}^d [x_k,x_k+h_k]$, set 
  \begin{equation*}
 \operatorname{Diff}_R(m)=\operatorname{Diff}(m,1,h_1,x)\cdots\operatorname{Diff}(m,d,h_d,x),\qquad x=(x_1,\ldots,x_d).
  \end{equation*}
 Define 
   \begin{equation*}
  \norm m .\operatorname{Var}_q(Q).=\sup_{\mathcal P} \Bigl[ \sum_{R\in\mathcal P} \abs{\operatorname{Diff}_R(M)}^q \Bigr]^{1/q}, \quad 0<q<\infty.
   \end{equation*}
  The supremum is formed over all partitions $\mathcal P$ of the rectangle $Q$ into subrectangles.
  
  Given $1\le{}k<{}d$, and  $y=(y_1,\ldots,y_k)\in \mathbb R^k$, and a map $\alpha\mid \{1,\ldots,k\}\to \{1,\ldots,d\}$, let 
  $m_{y,\alpha}$ be the function from $\mathbb R^{d-k}$ to $\mathbb C$ obtained from $m$ by restricting the $\alpha(j)$th coordinate to be $y_j$, $1\le j\le k$. Then, the $V_q(Q)$ norm is 
  \begin{equation*}
 \norm m.V_q(Q).=\norm m.\infty.+\norm m.\operatorname{Var}_q(\mathbb R^d).+\sup_{k,\alpha}\sup_{y\in\mathbb R^k}
 \norm m_{k,\alpha}.\operatorname{Var}_q (Q_{y,\alpha}). .
  \end{equation*}
 Here, we let $Q_{y,\alpha}$ be the cube obtained from $Q$ by restricting the $\alpha(j)$th coordinate to be $y_j$, $1\le{}j\le{}k$. 
 
  Recall the notation $\Delta^d$ for the lacunary intervals in $d$ dimensions, and in particular (\ref{e.higherlp}).  
 
   \begin{theorem}\label{t.crsIndim}  Suppose that $\abs{\frac12-\frac1p}<\frac1q$.  For functions $m\mid \mathbb R^d\longrightarrow \mathbb C$, we have the estimate on the multiplier norm of $m$  
   \begin{equation*}
  \norm \operatorname A_m.p.\lesssim{}\sup_{R\in\Delta^d} \norm m\ind R.V_q(\mathbb R^d).
   \end{equation*}
   \end{theorem} 
  
  \subsection{Proof of \thm.crsIndim/}
  
  The argument is a reprise of that for the one dimensional case.   We begin with definitions in one dimension. 
  Let $B$ be a linear space with norm $\norm \, .B.$.  For an interval $I$ let $\mathcal E(I,B)$ be the linear space of step functions 
  $m\mid I\longrightarrow B$ with finite range.  Thus, 
   \begin{equation*}
  m=\sum_{j=1}^J b_j \ind {I_j}
   \end{equation*}
  for a  finite partition $\{I_1,\ldots,I_j\}$ of $I$ into intervals, and a sequence of values $b_j\in B$.  If $B=\mathbb C$, we write simply 
  $\mathcal E(I)$.  There are a family of norms that we impose on $\mathcal E(I,B)$.  
   \begin{equation*}
  \langle\!\langle m\rangle\!\rangle_{\mathcal E(I,B),q}=\Bigl[ \sum_{j=1}^J \lVert b_j \rVert_B^q \Bigr]^{1/q},\qquad 1\le{}q\le\infty.
   \end{equation*}
    
    For a rectangle $R=R_1\times\cdots R_d$, set 
     \begin{equation*}
    \mathcal E(R):=\mathcal E(R_1, \mathcal E(R_2,\cdots ,\mathcal E(R_d,\mathbb C)\cdots))
     \end{equation*}
    The following Lemma is a variant of \l.decouple/. 
    
     \begin{lemma}\label{l.DEcouple}  Let $m\mid \mathbb R^d\longrightarrow\mathbb C$ be a function such that $m\ind R\in\mathcal E(R)$ for all $R\in\Delta^d$. Then, 
   we have these two estimates for the multiplier $\operatorname A_m$. 
    \begin{gather} \label{e.DE1} 
   \norm \operatorname A_m .p.\lesssim\sup_{R\in\Delta^d}\langle\!\langle m\rangle\!\rangle_{\mathcal E(R),2},\qquad 1<p<\infty,
   \\
   \label{e.DE2} 
   \norm \operatorname A_m .p.\lesssim\sup_{R\in\Delta^d}\langle\!\langle m\rangle\!\rangle_{\mathcal E(I,B),q},\qquad 1<p<\infty,\ \abs{\tfrac12-\tfrac1p}<\tfrac1q.
    \end{gather}
    \end{lemma}

   \begin{proof} The first claim, the obvious bound at $L^2$, and complex interpolation prove the second claim. 
 
  As for the first claim, take  a multiplier $m$ for which the right hand side in (\ref{e.DE1}) is $1$. To each $R\in\Delta^d$, there is a partition $\Omega_R$ of $R$ into 
  a finite number of rectangles so that 
   \begin{gather*}
  m \ind R {}={}\sum_{\omega\in\Omega_R} a_\omega \ind \omega ,
  \\
  \sum_{\omega\in\Omega_R}\abs{a_\omega}^2\le1.
   \end{gather*}
  This conclusion is obvious for $d=1$, and induction on dimension will prove it in full generality.

  Then observe that by Cauchy-Schwarz,  
   \begin{equation*}
  \abs{\operatorname S_R \operatorname A_m}\le{}\operatorname S^{\Omega_R}.
   \end{equation*}
  Set $\Omega=\bigcup_{R\in\Delta^d}\Omega_R$.  Using the Littlewood Paley inequality (\ref{e.higherlp}), and Rubio de Francia's inequality in $d$ dimensions, we may estimate 
   \begin{align*}
  \norm \operatorname A_m f.p.\simeq{}& \norm \operatorname   S^{\Delta^d}\operatorname A_m f .p.
  \\{}\lesssim{}& \norm \operatorname S^\Omega f.p.
  \\{}\lesssim{}& \norm  f.p..
   \end{align*}
  The last step requires that $2\le{}p<\infty$, but the operator norm  $\norm \operatorname A_m.p.$ is invariant under conjugation of $p$, so that 
  we need only consider this range of $p$'s.
  
  \end{proof}
  
  We extend the notion of $\mathcal E(I,B)$.  Let $B$ be a Banach space, and set 
  $
  \mathcal U_q(I,B)
  $ to be the Banach space of functions $m\mid I\longrightarrow B$ for which the norm below is finite. 
   \begin{equation*}
  \norm m .\mathcal U_q.:=\inf\Bigl\{ \sum_j \langle\!\langle m_j\rangle\!\rangle_{\mathcal E(I,B),q}\mid m=\sum_j m_j,\ m_j\in\mathcal E(I,B)\Bigr\}.
  \end{equation*}
 For a rectangle $R=R_1\times\cdots R_d$, set 
   \begin{equation*}
    \mathcal U_q(R):=\mathcal U_q(R_1, \mathcal U_q(R_2,\cdots ,\mathcal U_q(R_d,\mathbb C)\cdots)).
     \end{equation*}
  By convexity, we clearly have the inequalities 
   \begin{gather*}
  \norm \operatorname A_m.p.\lesssim\sup_{R\in\Delta^d}\norm m\ind R.\mathcal U_2(R).,\qquad 1<p<\infty,
  \\
  \norm \operatorname A_m .p.\lesssim\sup_{R\in\Delta^d}\norm m\ind R.\mathcal U_q(R)., \qquad 1<p<\infty,\ \abs{\tfrac12-\tfrac1p}<\tfrac1q.
   \end{gather*}
  As well, we have the inclusion $\mathcal U_q(R)\subset\operatorname{Var}_q(R)$, for $1\le{}q<\infty$.  The reverse inclusion is not 
  true in general, nevertheless the inclusion is true with a small perturbation of indicies.  
  
  Let us note that the definition of the $q$ variation space on an interval, given in (\ref{e.VAR}), has an immediate extension to a setting in which the functions $m$ take values in a Banach space $B$.  Let us denote this space as $\operatorname{V}_q(I, B)$. 
  
   \begin{lemma}\label{l.inclusion}  For all $1\le{}p<q<\infty$, all intervals $I$, and Banach spaces $B$, we have the inclusion 
   \begin{equation*}
  \operatorname{V}_p(I,B)\subset \mathcal U_q(I,B).
   \end{equation*}
  For all pairs of intervals $I,J$, we have 
   \begin{equation*}
  V_q(I\times J, B)=V_q(I, V_q(J,B)).
   \end{equation*}
  In addition,  for all  rectangles $R$, we have 
   \begin{equation*}
  \operatorname{V}_p(R)\subset\mathcal U_q(R).
   \end{equation*}
  In each instance, the inclusion map is bounded. 
  \end{lemma} 
 
 The first claim of the Lemma is proved by a trivial modification of the proof of \l.Vq/. (The martingale convergence theorem 
 holds for all Banach space valued martingales.) The second claim is  easy to verify, and the last claim is a corollary to the first two.

 \section{Notes and Remarks } 
 
  \begin{remark}\label{r.carleson}  L.~Carleson \cite{carleson} first noted the possible extension of the Littlewood Paley inequality,  proving in 1967 that \thm.lp/ holds in  the special case that $\Omega=\{[n,n+1)\mid n\in\mathbb Z\}$. He also noted that the inequality does not extend  to $1<p<2$. 
 A corresponding 
 extension to homethetic parallelepipeds was given by A.~C{\'o}rdoba \cite{MR83i:42015}, who also pointed out the connection to multipliers.
  \end{remark}
 
  \begin{remark}\label{r.approach}   Rubio de Francia's paper \cite{rubio} adopted an approach that we could outline this way.  The reduction to the well distributed case is made, and we have borrowed that line of reasoning from him.  This permits one  to define a smooth operator $\operatorname G^\Omega$ in (\ref{e.G}).  That 
$\operatorname G^\Omega$ is bounded on $L^p$, for $2<p<\infty$, is a consequence of a bound on the sharp function.  In our notation, that sharp function estimate would be 
 \begin{equation} \label{e.sharp} 
(\operatorname G^\Omega f)^\sharp\lesssim{} (\operatorname  M\abs{f}^2)^{1/2}.
 \end{equation}
The sharp   denotes the   function 
 \begin{equation*} 
g^\sharp=\sup_{J} \ind J  \Bigl[ \dashint_J \Abs{ g-\dashint_J g}^2 \; dx \Bigr]^{1/2},
 \end{equation*}
the supremum being formed over all intervals $J$.
It is known that $\norm g.p.\simeq \norm g^\sharp .p.$ for $1<p<\infty$.  Notice that our proof can be  adapted to prove a dyadic version of (\ref{e.sharp}) for the tile operator  $\operatorname T^\Omega$ if desired. The sharp function estimate has the advantage of quickly giving weighted inequalities. It has the disadvantage of not easily generalizing to higher dimensions.  On this point, see  R.~Fefferman \cite{MR90e:42030}.
 \end{remark}

 \begin{remark}\label{r.1weighted}  The weighted version of Rubio de Francia's inequality states that for all weights $w\in A^1$, one has 
 \begin{equation*}
\int \abs{\operatorname S^\Omega f}^2\; w\; dx\lesssim{}\int \abs{f}^2 \; w\; dx.
 \end{equation*}
There is a similar conclusion for multipliers. 
 \begin{equation*}
\int \abs{\operatorname A_m f}^2\; w\; dx\lesssim{}\sup_{R\in\Delta^d} \lVert m\ind R \rVert_{\operatorname{V}_2(R)}^2\int \abs{f}^2 \; w\; dx.
 \end{equation*}
See Coifman, Rubio de Francia, and Semmes \cite{MR89e:42009}, for one dimension and Q.~Xu \cite{MR98b:42024} for  dimensions greater than $1$. 
 \end{remark}

 \begin{remark}\label{r.alternate} Among those authors who made a contribution to Rubio de Francia's one dimensional inequality, 
  P.~Sj{\''o}lin \cite{MR88e:42018} provided an alternate derivation of Rubio de Francia's sharp function estimate (\ref{e.sharp}).  
  In another direction, observe that (\ref{e.duality}), Rubio de Francia's inequality has the dual formulation $\norm f.p.\lesssim{}\norm \operatorname S^\Omega f.p.$, 
  for $1<p<2$.  J.~Bourgain \cite{MR87m:42008} established a dual endpoint estimate for the unit circle:
   \begin{equation*}
  \norm f.H^1.\lesssim{}\norm \operatorname S^\Omega f.1.
   \end{equation*}
  This inequality in higher dimensions seems to be open. 
   \end{remark} 
  
   \begin{remark}\label{r.p<2}  Rubio de Francia's inequality does not extend below $L^2$.  While this is suggested by the duality estimates, an explicit example is given in one dimension  by $\widehat f=\ind {[0,N]}$, for a large integer $N$, and $\Omega=\{(n,n+1)\mid n\in\mathbb Z\}$.  Then, it is easy to see that 
   \begin{equation*}
  N^{1/2}\ind {[0,1]}\lesssim{}\operatorname S^\Omega f,\qquad \norm f.p.\simeq N^{p/(p-1)}, \quad 1<p<\infty.
   \end{equation*}
  It follows that $1<p<2$ is not permitted in Rubio de Francia's inequality. 
   \end{remark} 
  
   \begin{remark}\label{r.p<<2}  In considering an estimate below $L^2$,
  in any dimension, we have the following interpolation argument available to us for all well distributed collections $\Omega$. We have the estimate 
   \begin{equation} \label{e.max}  
  \sup_{s\in{\mathbf T}^\Omega} \ind {R_s}
  \frac{ \abs{\ip f, \varphi_s .}}{\sqrt{\abs{R_s}}} \lesssim{}\operatorname  M f
   \end{equation}
  where $M$ denotes the strong maximal function.  Thus, the right hand side is a bounded operator on all $L^p$.  By taking a value of $p$ very close to one, and interpolating with the $L^2$ bound for $\operatorname S^\Omega$, we see that 
    \begin{equation} \label{e.p<2}
   \biggl\lVert   \NOrm \ind {R_s} \frac{ \abs{\ip f, \varphi_s .}}{\sqrt{\abs{R_s}}} . \ell^q(\Omega) . \biggr\rVert_p\lesssim{}\norm f.p.,\qquad 1<p<2,\ 
   \frac p{p-1}<q<\infty.
    \end{equation}
  By (\ref{e.lpvector}), this inequality continues to hold without the well distributed assumption. 
  Namely, using (\ref{e.vectorH}), for all disjoint collections of rectangles $\Omega$, 
   \begin{equation*}
  \lVert \norm \operatorname S_\omega f.\ell^q(\Omega). \rVert_p\lesssim\norm f.p.\qquad 1<p<2,\ \frac p{p-1}<q<\infty. 
   \end{equation*}
   \end{remark}

  \begin{remark}\label{r.p3} 
  Cowling and Tao \cite{MR2002j:42028}, for $ 1<p<2$, construct $ f\in L^p$ 
  for which 
\begin{equation*}
\norm  { \norm \operatorname S _{\omega }f.\ell^{p'}(\Omega ). } .L ^{p}.
=\infty \,, 
\end{equation*}
 forbidding one possible extension of the interpolation above to a natural endpoint 
estimate.  
\end{remark}

\begin{remark}\label{r.quek}   Quek \cite{MR1950722}, on the other hand, 
finds a sharp endpoint estimate
\begin{equation*}
\norm  { \norm \operatorname S _{\omega }f.\ell^{p'}(\Omega ). } .L ^{p,p'}.
\lesssim \norm f.p.,\qquad 1<p<2,\ p'=\tfrac p{p-1}.
\end{equation*}
In this last inequality, the space $ L^{p,p'}$ denotes a Lorenz space.  
Indeed, he uses the weak $ L^1$ estimate (\ref{e.max}), together with a complex 
interpolation method. 
\end{remark}

   \begin{remark}\label{r.higher} The higher dimensional formulation of Rubio de Francia's inequality did not admit an immediately clear formulation.  J.-L.~Journ\'e 
  \cite{MR88d:42028} established the Theorem in the higher dimensional case, but used a very sophisticated proof.  Simpler arguments, very close in spirit to what we have presented, were given by  F.~Soria \cite{MR88g:42026} in two dimensions, and in higher dimensions  by S.~Sato \cite{MR92c:42020} and X.~Zhu \cite{MR93f:42041}. These authors continued to focus on the $G$ function (\ref{e.G}), instead of the time frequency approach we have used. 
   \end{remark} 
  
   \begin{remark}\label{r.simplehigher}  We should mention that if one is considering the higher dimensional version of \thm.lp/, with the simplification that  the collection of rectangles consists only of cubes, then the method of proof need not invoke the difficulties of the $BMO$ theory of Chang and Fefferman. The usual one parameter $BMO$ theory will suffice.  The same comment holds if all the rectangles in $\Omega$ are homeothetic under translations and application of a power of a fixed expanding matrix. 
   \end{remark}

   \begin{remark}\label{r.other}   It would be of interest to establish variants of Rubio de Francia's inequality for other collections of sets in the plane. 
  A.~Cordoba has established a preliminary result in this direction for finite numbers of sectors in the plane.     
  G.~Karagulyan and the author \cite {math.CA/0404028}  provide a more  result for more general sets of directions, 
  presuming \emph{a priori} bounds on the maximal function associated with this set of 
  directions. 
  \end{remark}

   \begin{remark}\label{r.weighted} The  inequality (\ref{e.vectorH}) is now typically seen as a consequence of the general theory of weighted inequalities. 
  In particular, if $h\in L^1(\mathbb R)$, and $\epsilon>0$, it is the case that $ (M h)^{1-\epsilon}$ is a weight in the Mockenhoupt class $A^1$. In particular, this observation implies (\ref{e.vectorH}). See the material on weighted inequalities in 
  E.M.~Stein \cite{stein}.
   \end{remark}

   \begin{remark}\label{r.jose}  Critical to the proof of Rubio de Francia's inequality is the $L^2$ boundedness 
  of the tile operator $\operatorname T^\Omega$. This is of course an immediate consequence of the well distributed assumption.  It would be of 
  some interest to establish a reasonable geometric condition on the tiles which would be sufficient for the $L^2$ boundedness of the 
  operator $\operatorname T^\Omega$.  In this regard, one should consult the inequality of J.~Barrionuevo and the author \cite{MR1955263}. 
  This inequality    is of a weak type, but is sharp. 
   \end{remark}

   \begin{remark}\label{r.olev}  V.~Olevskii \cite{MR95a:42012a} independently established a version of \thm.crs/ on the integers. 
   \end{remark}

  \begin{remark}\label{r.multipliers}  Observe that, in a certain sense, the multiplier result \thm.crs/ is optimal.  In one dimension,  let $\psi$ denote  a smooth 
bump function $\psi$ with frequency support in $[-1/2,1/2]$, and for random choices of signs $\varepsilon_k$, $k\ge1$, and integer $N$, consider the multiplier 
 \begin{equation*}
m=\sum_{k=1}^N \varepsilon_k \Tr k \widehat\psi 
 \end{equation*}
Apply this multiplier to the function $\widehat f=\ind {[0,N]}$.  By an application of Khintchine's inequality, 
 \begin{equation*}
\mathbb E \norm \operatorname A_m f.p.\simeq \sqrt N, \qquad 1<p<\infty.
 \end{equation*}
On the other hand, it is straight forward to verify that $\norm f.p.\simeq N^{p/(p-1)}$.  We conclude that 
 \begin{equation*}
\norm \operatorname A_m.p.\gtrsim N^{\abs{1/2-1/p}}.
 \end{equation*}
Clearly $\mathbb E\norm m.V_q.\simeq N^{1/q}$.  That is, up to arbitrarily small constant, the values of $q$ permitted in \thm.crs/ are optimal. 
 \end{remark}

 \begin{remark}\label{r.olevskii} 
 V.~Olevskii \cite{MR98e:42007} refines the notion in which \thm.crs/ is  optimal.  The argument is phrased in terms of multipliers for $\ell^p(\mathbb Z)$.   
 It is evident  that the $q$ variation norm is preserved under  homeomorphims.   That is   let $\phi\mid \mathbb T\longrightarrow \mathbb T$ be a homeomorphism. Then    $\norm m.V_q(\mathbb T).=\norm m\circ \phi .V_q(\mathbb T).$.  For a multiplier $m\mid \mathbb T\longrightarrow \mathbb R$, let 
 \begin{equation*}
 \norm m .M^0_p.=\sup_\phi \sup_{\norm f.\ell^p(\mathbb Z).=1} \NOrm  \int_\mathbb T \widehat f (\tau) m\circ\phi(\tau)e^{in\tau}\;d\tau .p.  
 \end{equation*}
Thus, $M^0_p$ is the supremum over multiplier norms of $m\circ \phi$, for all homeomorphims $\phi$.   Then, Olevskii shows that 
if $\norm m.M^0_p.<\infty$, then $ m$ has finite $q$ variation for all $\abs{\frac12-\frac1p}<\frac1q$. 
 \end{remark}

   \begin{remark}\label{r.bg}  E.~Berkson and T.~Gillespie \cites{bg1,bg2,bg3} have extended the Coifman, Rubio, Semmes result to a setting 
  in which one has an operator with an appropriate spectral representation.  
   \end{remark}
  
   \begin{remark}\label{r.hare}   The Rubio de Francia inequalities are in only one direction.  K.E.~Hare and I.~Klemes \cites{hareklemes1,hareklemes2,hareklemes3} have undertaken a somewhat general study of necessary 
  and sufficient conditions on a class of intervals to satisfy a the inequality that is reverse to that of 
  Rubio.  A theorem from \cite{hareklemes3} concerns a collection of intervals $\Omega=\{\omega_j\mid j\in\mathbb Z\}$ which are assumed to partition $\mathbb R$, and satisfy $\abs {\omega_{j+1}}/\abs{ \omega_j}\to\infty$.  Then one has 
   \begin{equation*}
  \norm f.p.\lesssim\norm \operatorname S^\Omega f.p.,\qquad 2<p<\infty.
   \end{equation*}
  What is important here is that the locations of the $\omega_j$ are not specified.  
  The authors conjecture that it is enough to have $\abs {\omega_{j+1}}/\abs{ \omega_j}>\lambda>1$.
   \end{remark}

\begin{remark}\label{r.earl} Certain operator theoretic variants of issues related to 
Rubio de Francia's inequality  are discussed in the papers of Berkson and Gillespie \cites{MR2083875
,
MR2054306
,
MR2001939}.
\end{remark}

\end{document}